

\magnification=1100
\overfullrule0pt

\input prepictex
\input pictex
\input postpictex
\input amssym.def


\def\qed{\hbox{\hskip 1pt\vrule width4pt height 6pt depth1.5pt \hskip 1pt}}

\def\CC{{\Bbb C}}
\def\RR{{\Bbb R}}
\def\ZZ{{\Bbb Z}}

\def\cB{{\cal B}}

\def\cF{{\cal F}}

\def\Card{{\rm Card}}

\def\Ind{{\rm Ind}}
\def \supp{{\rm supp}}

 
\font\smallcaps=cmcsc10
\font\titlefont=cmr10 scaled \magstep1

\font\tinyrm=cmr10 at 8pt


\newcount\sectno
\newcount\subsectno
\newcount\resultno

\def\section #1. #2\par{
\sectno=#1
\resultno=0
\bigskip\noindent\centerline{\smallcaps #1.  #2}~\medbreak}

\def\subsection #1\par{\medskip\noindent{\bf  #1} }


\def\prop{ \global\advance\resultno by 1
\bigskip\noindent{\bf Proposition \the\sectno.\the\resultno. }\sl}
\def\lemma{ \global\advance\resultno by 1
\bigskip\noindent{\bf Lemma \the\sectno.\the\resultno. }
\sl}
\def\remark{ \global\advance\resultno by 1
\bigskip\noindent{\bf Remark \the\sectno.\the\resultno. }}
\def\example{ \global\advance\resultno by 1
\bigskip\noindent{\bf Example \the\sectno.\the\resultno. }}
\def\cor{ \global\advance\resultno by 1
\bigskip\noindent{\bf Corollary \the\sectno.\the\resultno. }\sl}
\def\thm{ \global\advance\resultno by 1
\bigskip\noindent{\bf Theorem \the\sectno.\the\resultno. }\sl}
\def\defn{ \global\advance\resultno by 1
\bigskip\noindent{\it Definition \the\sectno.\the\resultno. }\slrm}
\def\endthm{\rm\bigskip}

\def\endprop{\rm\bigskip}

\def\pf{\rm\bigskip\noindent{\it Proof. }}
\def\endpf{\qed\hfil\bigskip}


\def\formula{\global\advance\resultno by 1
\eqno{(\the\sectno.\the\resultno)}}
\def\formulano{\global\advance\resultno by 1 (\the\sectno.\the\resultno)}
\def\tableno{\global\advance\resultno by 1
\the\sectno.\the\resultno. }
\def\lformula{\global\advance\resultno by 1
\leqno(\the\sectno.\the\resultno)}

\def\monthname {\ifcase\month\or January\or February\or March\or April\or
May\or June\or
July\or August\or September\or October\or November\or December\fi}

\newcount\mins  \newcount\hours  \hours=\time \mins=\time
\def\now{\divide\hours by60 \multiply\hours by60 \advance\mins by-\hours
     \divide\hours by60         
     \ifnum\hours>12 \advance\hours by-12
       \number\hours:\ifnum\mins<10 0\fi\number\mins\ P.M.\else
       \number\hours:\ifnum\mins<10 0\fi\number\mins\ A.M.\fi}


\nopagenumbers
\def\runningtitle{\smallcaps Skew shape representations}
\headline={\ifnum\pageno>1\eoheadline\else\firstheadline\fi}
\def\names{\smallcaps arun ram}
\def\firstheadline{\noindent Preliminary Draft \hfill  \today}
\def\firstheadline{}
\def\eoheadline{\ifodd\pageno\oddheadline\else\evenheadline\fi}
\def\oddheadline{\tenrm\hfil\runningtitle\hfil\folio}
\def\evenheadline{\tenrm \folio\hfil{\names}\hfil}

\vphantom{$ $}  
\vskip.75truein

\centerline{\titlefont Skew shape representations are irreducible}
\bigskip
\centerline{\rm Arun Ram${}^\ast$}
\centerline{Department of Mathematics}
\centerline{Princeton University}
\centerline{Princeton, NJ 08544}
\centerline{{\tt rama@math.princeton.edu}}
\centerline{Preprint:  August 4, 1998}

\footnote{}{\tinyrm ${}^\ast$ Research supported in part by National
Science Foundation grant DMS-9622985, and a Postdoctoral Fellowship
at Mathematical Sciences Research Institute.}

\bigskip

\noindent{\bf Abstract.}
In this paper all of the classical constructions of A.~Young 
are generalized to affine Hecke algebras of type A.
It is proved that the calibrated irreducible representations of the affine 
Hecke algebra are indexed by placed skew shapes
and that these representations can be constructed explicitly with a 
generalization of Young's seminormal construction of the irreducible
representations of the symmetric group.  The seminormal construction
of an irreducible calibrated module does not produce a basis on which 
the affine Hecke algebra acts integrally
but using it one is able to pick out a different basis,
an analogue of Young's natural basis,  which does generate an integral lattice
in the module.   Analogues of the ``Garnir relations'' play an important
role in the proof.  The Littlewood-Richardson coefficients arise naturally
as the decomposition multiplicities for the restriction of an irreducible
representation of the affine Hecke algebra to the Iwahori-Hecke algebra.

\section 0. Introduction

My recent work [Ra3], [Ra4] on the representations of affine Hecke algebras has
been strongly motivated by the classical theory of Young tableaux.  
This research has resulted in the generalization
of many of A. Young's constructions to general finite root systems.
With these generalizations of standard Young tableaux 
one is able to use Young's classical ``seminormal construction''
to construct irreducible representations of affine Hecke algebras
corresponding to arbitrary finite crystallographic root systems.

Because the classical combinatorics of Young tableaux is much more advanced
than that of the newly developed generalization it is often possible
to give simpler proofs and more extensive results for the case of affine Hecke
algebras of type A.  The purpose of this paper is to compile some of these
results and proofs.  In particular, we obtain a generalization
of Young's natural basis and derive certain induction and restriction
rules which are not yet available in the general case.
It will also be more clear from the exposition here, how the 
generalization of the Young tableau theory given in [Ra3] and [Ra4]
relates to the classical setup, something which is not always
obvious when working in the general root system context.

\medskip\noindent
{\it The main results}
\medskip

\item{(1)} {\sl The definition of calibrated representations, and the
classification and construction of all irreducible calibrated representations
of the affine Hecke algebra type A. } 

\smallskip\noindent
These representations are indexed
by placed skew shapes.  The dimension of an irreducible calibrated
representation is the number of standard tableaux of the corresponding
skew shape and the representation is constructed by explicit
formulas which give the action of each generator of the affine Hecke
algebra on a specific basis, the elements of which are indexed by
standard Young tableaux.
This is a generalization of the constructions of A. Young [Y],
P. Hoefsmit [Ho], H. Wenzl [Wz], and Ariki and Koike [AK] (see [Ra2] for a
review of some of the unpublished results of Hoefsmit). 
Parts of Theorem 4.1 were first discovered by Cherednik and are stated (without
proof) in [Ch].  I am grateful to A. Zelevinsky for pointing this out to me and
to I. Cherednik for some informative  discussions.

\smallskip
\item{(2)}  {\sl The definition of an analogue of Young's natural basis
for each irreducible calibrated representation.}

\smallskip\noindent
Young's natural basis is the one that is most often used in the
study of irreducible representations of the symmetric group, it is
the one that is usually taken as the basis of the ``Specht module'',
see for example [JK], [Sg], [Fu].  It has the wonderful
property that it is an {\it integral} basis for the module,
i.e.  the matrices representating the action of the symmetric group
on this basis contain integer entries.  This is especially important
because it opens the door to a combinatorial study of the modular
representations of the symmetric group.

Using the analogue of the seminormal basis for the irreducible
calibrated representations of the affine Hecke algebra we can
define an analogue of Young's natural basis in each of these
representations.  As desired, this basis is an integral
basis of the module; the matrices representing the action 
of the affine Hecke algebra on this basis have all entries in the
ring $\ZZ[q,q^{-1}]$.  These results are a $q$-analogue of some of the
results in [GW].

One of the pleasant surprises one has when generalizing
Young's natural basis from this point of view is that
the ``Garnir relations'' take a particularly simple form:
If $\{v_L\}$ is Young's seminormal basis and $\{n_L\}$ is Young's natural basis
then the relations
$$v_L=0, \qquad\hbox{when $L$ is not a standard tableau,}$$
are the Garnir relations.  One recovers the Garnir relations in their
classical form by expanding the $v_L$ in terms of the $n_L$.

\smallskip
\item{(3)}  {\sl The classical Littlewood-Richardson coefficients
describe the decomposition of the restriction of an irreducible
representation of the affine Hecke algebra to the Iwahori-Hecke algebra.}

\smallskip\noindent
This result gives a completely new (and unexpected)
representation theoretic interpretation of the
Littlewood-Richardson coefficients.

\smallskip
\item{(4)}  {\sl Skew shapes arise naturally as indexes for the 
{\it irreducible} calibrated representations of the affine Hecke algebra
 of type A.}

\smallskip\noindent
Until now skew shapes have appeared in the combinatorial literature as 
something of a novelty, a useful combinatorial tool which indexes some
strangely well-behaved representations of the symmetric group.  
It has always been a surprise that the combinatorics
of the irreducible representations of the symmetric group generalizes
so beautifully to this special class of highly reducible representations
of the symmetric group.

This fact is no longer strange.  In fact, these representations
are {\it irreducible} representations of the affine Hecke algebra,
and thus are basic and fundamental.
Several of the skew Schur function identities in [Mac] I can be given
representation theoretic interpretations in this context, see Theorem 6.2 and
Corollary 6.3.

\medskip\noindent
{\it Acknowledgements}

This paper is part of a series [Ra3-5] [RR1-2] of papers on
representations of affine Hecke algebras.  During this work I have benefited from
conversations with many people. To choose only a few, there were discussions with
S. Fomin, F. Knop, L. Solomon, M. Vazirani and N. Wallach which played an
important role in my progress.  There were several times when
I tapped into J. Stembridge's fountain of useful knowledge about root systems. 
G. Benkart was a very patient listener on many occasions.  H. Barcelo, P. Deligne,
T. Halverson, R. Macpherson and R. Simion all gave large amounts of time to
let me tell them my story and every one of these sessions was helpful to me in
solidifying my understanding.

I single out Jacqui Ramagge with special thanks for everything she
has done to help with this project: from the most mundane typing and picture
drawing to deep intense mathematical conversations which helped to sort out
many pieces of this theory.  Her immense contribution is evident in
that some of the papers in this series on representations of affine Hecke
algebras are joint papers.  

A portion of this research was done during a semester stay at Mathematical
Sciences Research Institute where I was supported by a Postdoctoral
Fellowship.  I thank MSRI and National Science Foundation for support of
my research.

\bigskip

\section 1. The affine Hecke algebra of type A

\subsection{Affine braids.}  

There are three common ways of
depicting affine braids [Cr], [GL], [Jo]:
\smallskip
\item{(a)}  As braids in a (slightly thickened) cylinder,
\smallskip
\item{(b)}  As braids in a (slightly thickened) annulus,
\smallskip
\item{(c)}  As braids with a flagpole.

\smallskip\noindent
See Figure 1.
The multiplication is by placing one cylinder on
top of another,  placing one annulus inside another,
or placing one flagpole braid on top of another.
These are equivalent formulations:
an annulus can be made into a cylinder by turning up the edges,
and a cylindrical braid can be made into a flagpole braid
by putting a flagpole down the middle of the cylinder and 
pushing the pole over to the left so that the strings
begin and end to its right.

The group formed by the affine braids with $n$ strands is the
{\it affine braid group} 
$\tilde\cB_n$ of type A.  Let $\omega$, $T_i$ for
$0\le i\le n-1$, and $x_i$ for $1\le i\le n$, be as given in
Figure 2.
The following identities can be checked by drawing pictures:
$$\matrix{
\hbox{(a)\ \   $T_iT_j=T_jT_i$,}\hfill
&\qquad&\hbox{for $|i-j|>1$,}\hfill \cr
\cr
\hbox{(b)\ \   $T_iT_{i+1}T_i = T_{i+1}T_iT_{i+1}$,}\hfill
&&\hbox{for $0\le i\le n-1$,}\hfill\cr
\cr
\hbox{(c)\ \   $\omega T_i \omega^{-1} = T_{i-1}$,} \hfill
&&\hbox{for $0\le i\le n-1$,}\hfill\cr
\cr
\hbox{(d)\ \   $x_iT_j= T_jx_i$,} \hfill
&&\hbox{if $|i-j|>1$,}\hfill\cr
\cr
\hbox{(e)\ \   $x_{i+1}= T_ix_iT_i$,} \hfill
&&\hbox{for $1\le i\le n-1$,}\hfill\cr
\cr
\hbox{(f)\ \   $x_ix_j=x_jx_i$, } \hfill
&&\hbox{for $1\le i,j,\le n$,}\hfill \cr
\cr
\hbox{(g)\ \   $x_nx_1^{-1} = T_0T_{n-1}\cdots T_2T_1T_2\cdots T_{n-1}$,}
\hfill\cr
\cr
\hbox{(h)\ \   $x_n = \omega T_1T_2\cdots T_{n-1}$,} \hfill \cr
\cr
\hbox{(i)\ \   $\omega^n = x_1x_2\cdots x_n$,} \hfill\cr
}
\formula$$
where the indices on the elements $T_i$ are taken modulo $n$.
The elements $T_i$, $0\le i\le n-1$, and $\omega$ generate $\tilde \cB_n$.
The {\it braid group} is the subgroup $\cB_n$ generated by the $T_i$, $1\le i\le
n-1$. The elements $x_i$, $1\le i\le n$, generate an abelian group
$X\subseteq \tilde \cB_n$. If
$\gamma = (\gamma_1,\gamma_2,\ldots,\gamma_n)\in
\ZZ^n$ define
$$x^\gamma = x_1^{\gamma_1}x_2^{\gamma_2}\cdots x_n^{\gamma_n}.\formula$$
The symmetric group $S_n$ acts on $\ZZ^n$ by permuting the
coordinates.  This action induces an action on $X$ by
$$wx^\gamma = x^{w\gamma}, \qquad \hbox{for $w\in S_n$, $\gamma\in \ZZ^n$.}
$$

\subsection{The affine Hecke algebra.}

Fix an element $q\in \CC^*$ which is not a root of unity.  The {\it affine
Hecke algebra} $\tilde H_n$ is the quotient of the group algebra
$\CC\tilde\cB_n$ by the relations
$$T_i^2=(q-q^{-1})T_i+1,
\qquad 0\le i\le n.\formula$$
The images of $T_i$, $x_i$ and $\omega$ in $\tilde H_n$ are again denoted by
$T_i$, $x_i$ and $\omega$.  The Laurent polynomial ring $\CC[X] = \CC[x_1^{\pm 1},\ldots, x_n^{\pm 1}]$ is a
(large) commutative subalgebra of $\tilde H_n$.

The relations
$T_i^{-1} = T_i - (q-q^{-1})$ and $x_{i+1}=T_ix_iT_i$ can be used
to derive the identities
$$
x_{i+1}T_i = T_ix_i+(q-q^{-1})x_{i+1}, \qquad\hbox{and}\qquad
x_iT_i = T_ix_{i+1}-(q-q^{-1})x_{i+1}. 
\formula$$
More generally, if $\gamma = (\gamma_1,\gamma_2,\ldots,\gamma_n)\in \ZZ^n$ then
$$x^\gamma T_i = T_i x^{s_i\gamma} + (q-q^{-1}) 
{(x^\gamma-x^{s_i\gamma})x_{i+1}\over x_{i+1}-x_i },\formula$$
where $s_i\in S_n$ is the simple transposition $(i,i+1)$.
The right hand term in this expression can always be written
as a Laurent polynomial in $x_1,\ldots, x_n$.  This important relation is due
to Bernstein, Zelevinsky and Lusztig [Lu].
The affine Hecke algebra
$\tilde H_n$ can be defined as the algebra generated by $T_i$, $1\le i\le n$,
and $x_i$, $1\le i\le n$ subject to the relations in (1.1a), (1.1b), (1.1f),
(1.3) and (1.5).

\subsection{The symmetric group.}

The {\it simple transpositions} are the elements
$s_i=(i,i+1)$, $1\le i\le n-1$, in $S_n$.
A {\it reduced word} for a permutation $w\in S_n$ is an expression
$w=s_{i_1}\cdots s_{i_p}$ of minimal length.  This minimal length is called the
{\it length} $\ell(w)$ of $w$.  The symmetric group $S_n$ is partially
ordered by the Bruhat-Chevalley order:
$v\le w$ if a reduced expression $s_{i_1}\cdots s_{i_p}$ for
$w$ has a subword $s_{i_{k_1}}\cdots s_{i_{k_\ell}}$, 
$1\le k_1<\cdots<k_\ell\le p$ which is equal to $v$ in $S_n$.

\subsection{The Iwahori-Hecke algebra.}

The {\it Iwahori-Hecke algebra} $H_n$ is the subalgebra of $\tilde H_n$
generated by the elements $T_i$, $1\le i\le n-1$.  For each $w\in S_n$ let
$$T_w = T_{i_1}\cdots T_{i_p},\formula$$
where $w=s_{i_1}\cdots s_{i_p}$ is a reduced word for $w$.  Since the $T_i$
satisfy the braid relations (1.1a,b), the element
$T_w$ is independent of the choice of the reduced word of $w$.  The
elements $T_w$, $w\in S_n$, are a basis of $H_n$ [Bou, IV \S 2 Ex. 23].

%
%
%
%
%
%

\vfill\eject

\section 2. Tableau combinatorics

\subsection{Skew shapes and standard tableaux.}

A partition $\lambda$ is a collection of $n$
boxes in a  corner.  We shall conform to the conventions in [Mac] and assume
that gravity goes up and to the left. 
$$
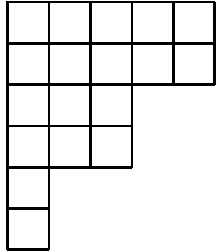
\beginpicture
\setcoordinatesystem units <0.5cm,0.5cm>         
\setplotarea x from 0 to 4, y from 0 to 3    
\linethickness=0.5pt                          
\putrule from 0 6 to 5 6          %
\putrule from 0 5 to 5 5          
\putrule from 0 4 to 5 4          %
\putrule from 0 3 to 3 3          %
\putrule from 0 2 to 3 2          %
\putrule from 0 1 to 1 1          %
\putrule from 0 0 to 1 0          %

\putrule from 0 0 to 0 6        %
\putrule from 1 0 to 1 6        %
\putrule from 2 2 to 2 6        %
\putrule from 3 2 to 3 6        
\putrule from 4 4 to 4 6        %
\putrule from 5 4 to 5 6        %
\endpicture
$$
Any partition $\lambda$ can be identified with
the sequence $\lambda=(\lambda_1\ge \lambda_2\ge \ldots )$
where $\lambda_i$ is the number of boxes in row $i$ of $\lambda$.
The rows and columns are numbered in the same way as for matrices.
In the example above we have $\lambda=(553311)$.
If $\lambda$ and $\mu$ are partitions such that $\mu_i\le \lambda_i$
for all $i$ we write $\mu\subseteq \lambda$.  
The {\it skew shape} $\lambda/\mu$
consists of all boxes of $\lambda$ which are not in $\mu$.
Any skew shape is a union of connected components.  Number the boxes of each
skew shape $\lambda/\mu$ along major diagonals from southwest to northeast and
$$\hbox{write ${\rm box}_i$ to indicate the box numbered $i$.}$$

Let $\lambda/\mu$ be a skew shape with $n$ boxes. 
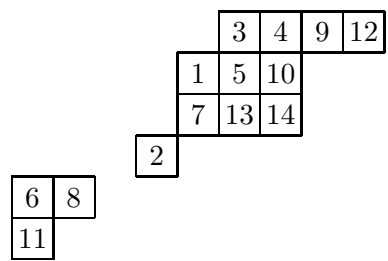
A {\it standard tableau of shape} $\lambda/\mu$ is a filling  
of the boxes in the skew shape
$\lambda/\mu$ with the numbers $1,\ldots,n$ such that 
the numbers increase from left to right in each row and from top to
bottom down each column. 
Let 
$$\cF^{\lambda/\mu}=\{\hbox{standard tableaux of shape
$\lambda/\mu$}\}.$$
The {\it column reading tableau} $C$ of shape $\lambda/\mu$ is the standard
tableau obtained by entering the numbers $1,2,\ldots,n$ consecutively 
down the columns of $\lambda/\mu$, beginning with the southwest most
connected component and filling the columns from left to right.
The {\it row reading tableau} $R$ of shape $\lambda/\mu$ is the standard
tableau obtained by entering the numbers $1,2,\ldots,n$ left to right
across the rows of $\lambda/\mu$, beginning with the northeast most
connected component and filling the rows from top to bottom.
In general, if $L$ is a standard tableau and $w\in S_n$ then $wL$ will denote the
filling of $\lambda/\mu$ obtained by permuting the entries of $L$ according to
the permutation $w$.

\prop {\rm [BW, Theorem 7.1]} 
Given a standard tableau $L$ of shape $\lambda/\mu$ define the 
{\it word} of $L$ to be permutation
$$
w_L=\left(\matrix{
1&\cdots & n \cr
L({\rm box}_1)&\ldots &L({\rm box}_n)
}\right)
$$ 
where $L({\rm box}_i)$ is the entry in ${\rm box}_i$ of $L$. 
Let $C$ and $R$ be the column reading and row reading tableaux
of shape $\lambda/\mu$, respectively.   
The map
$$\matrix{\cF^{\lambda/\mu} &\longrightarrow &S_n \cr
L &\longmapsto &w_L \cr
}$$
defines a bijection from $\cF^{\lambda/\mu}$ to the interval
$[w_C, w_R]$ in $S_n$ (in the Bruhat-Chevalley order).
\endprop

\subsection{Placed skew shapes.}

Let $\RR+i[0,2\pi/\ln(q^2))=\{ a+bi\ |\ a\in \RR, 0\le b\le 2\pi/\ln(q^2)
\}\subseteq \CC$.  If $q$ is a positive real number then
the function
$$\matrix{
\RR+i[0,2\pi/\ln(q^2)) &\longrightarrow &\CC^*\hfill \cr
x &\longmapsto &q^{2x}=e^{\ln(q^2)x} \cr}$$
is a bijection.
The elements of $[0,1)+i[0,2\pi/\ln(q^2))$ index the $\ZZ$-cosets in 
$\RR+i[0,2\pi/\ln(q^2))$.

A {\it placed skew shape} is a pair $(c,\lambda/\mu)$
consisting of a skew shape $\lambda/\mu$ and a {\it content
function}
$$c\colon \{\hbox{boxes of $\lambda/\mu$}\} \longrightarrow
\RR+i[0,2\pi/\ln(q^2))
\qquad\hbox{such that}$$
$$
\matrix{
c({\rm box}_j)\ge c({\rm box}_i), \hfill
&\quad &\hbox{if $i<j$ and $c({\rm box}_j)-c({\rm box}_i)\in \ZZ$}, \cr
c({\rm box}_j)=c({\rm box}_i)+1,\hfill
&&\hbox{if and only if ${\rm box}_i$ and ${\rm box}_j$
are on adjacent diagonals, and} \cr
c({\rm box}_i)=c({\rm box}_j),\hfill
&&\hbox{if and only if ${\rm box}_i$ and ${\rm box}_j$
are on the same diagonal.} \cr
}$$
This is a generalization of the usual notion of the content of a box
in a partition (see [Mac] I \S 1 Ex. 3).

Suppose that $(c,\lambda/\mu)$ is a placed skew shape such that
$c$ takes values in $\ZZ$.  One can visualize $(c,\lambda/\mu)$
by placing
$\lambda/\mu$ on a piece of infinite graph paper 
where the diagonals of the 
graph paper are indexed consecutively (with elements of $\ZZ$) from southeast 
to northwest.  The {\it content} of a box $b$ is the index $c(b)$ of the
diagonal that $b$ is on.
In the general case, when $c$ takes values in 
$\RR+i[0,2\pi/\ln(q^2))$, one imagines 
a book with $r$ pages of infinite graph paper
where the diagonals of the graph paper are indexed consecutively 
(with elements of $\ZZ$) from southeast to northwest.  
The pages are numbered by values $\beta_1,\ldots,\beta_r$ 
from the set $[0,1)+i[0,2\pi/\ln(q^2))$  
and there is a skew shape $\lambda^{(k)}/\mu^{(k)}$ placed on page $\beta_k$.   
The skew shape $\lambda/\mu$ is a union of the disjoint skew shapes
$\lambda^{(i)}/\mu^{(i)}$,
$$\lambda/\mu = \lambda^{(1)}/\mu^{(1)}\cup\cdots\cup
\lambda^{(r)}/\mu^{(r)},$$
and the content function is given by
$$
c(b) = 
\hbox{(page number of the page containing $b$)}
+ \hbox{(index of the diagonal containing $b$).} 
$$

\subsection{Example.}

The following diagrams illustrate standard tableaux and
the numbering of boxes in a skew shape $\lambda/\mu$.
$$
\matrix{
\beginpicture
\setcoordinatesystem units <0.5cm,0.5cm>         
\setplotarea x from 0 to 4, y from 0 to 3    
\linethickness=0.5pt                          
\put{1} at 0.5 0.5
\put{2} at 0.5 1.5
\put{3} at 1.5 1.5
\put{4}  at 3.5 2.5
\put{5} at 4.5 3.5
\put{6} at 4.5 4.5
\put{7}  at 5.5 3.5
\put{8}  at 5.5 4.5
\put{10}  at 5.5 5.5
\put{9}  at 6.5 3.5
\put{11}  at 6.5 4.5
\put{12}  at 6.5 5.5
\put{13}  at 7.5 5.5
\put{14}  at 8.5 5.5
\putrule from 5 6 to 9 6          %
\putrule from 4 5 to 9 5          
\putrule from 4 4 to 7 4          %
\putrule from 3 3 to 7 3          %
\putrule from 3 2 to 4 2          %
\putrule from 0 2 to 2 2          %
\putrule from 0 1 to 2 1          %
\putrule from 0 0 to 1 0          %
\putrule from 0 0 to 0 2        %
\putrule from 1 0 to 1 2        %
\putrule from 2 1 to 2 2        %
\putrule from 3 2 to 3 3        
\putrule from 4 2 to 4 5        %
\putrule from 5 3 to 5 6        %
\putrule from 6 3 to 6 6        %
\putrule from 7 3 to 7 6        %
\putrule from 8 5 to 8 6        %
\putrule from 9 5 to 9 6        %
\endpicture
&\qquad\qquad
&
\beginpicture
\setcoordinatesystem units <0.5cm,0.5cm>         
\setplotarea x from 0 to 4, y from 0 to 3    
\linethickness=0.5pt                          
\put{11} at 0.5 0.5
\put{6} at 0.5 1.5
\put{8} at 1.5 1.5
\put{2}  at 3.5 2.5
\put{7} at 4.5 3.5
\put{1} at 4.5 4.5
\put{13}  at 5.5 3.5
\put{5}  at 5.5 4.5
\put{3}  at 5.5 5.5
\put{14}  at 6.5 3.5
\put{10}  at 6.5 4.5
\put{4}  at 6.5 5.5
\put{9}  at 7.5 5.5
\put{12}  at 8.5 5.5
\putrule from 5 6 to 9 6          %
\putrule from 4 5 to 9 5          
\putrule from 4 4 to 7 4          %
\putrule from 3 3 to 7 3          %
\putrule from 3 2 to 4 2          %
\putrule from 0 2 to 2 2          %
\putrule from 0 1 to 2 1          %
\putrule from 0 0 to 1 0          %
\putrule from 0 0 to 0 2        %
\putrule from 1 0 to 1 2        %
\putrule from 2 1 to 2 2        %
\putrule from 3 2 to 3 3        
\putrule from 4 2 to 4 5        %
\putrule from 5 3 to 5 6        %
\putrule from 6 3 to 6 6        %
\putrule from 7 3 to 7 6        %
\putrule from 8 5 to 8 6        %
\putrule from 9 5 to 9 6        %
\endpicture
\cr
\hbox{$\lambda/\mu$ with boxes numbered}
&&\hbox{A standard tableau $L$ of shape $\lambda/\mu$} \cr
}
$$
The word of the standard tableau $L$ is the permutation
$w_L=(11,6,8,2,7,1,13,5,14,3,10,4,9,12)$ (in one-line 
notation).

The following picture shows the contents of the boxes in the
placed skew shape $(c,\lambda/\mu)$ such that the sequence
$(c({\rm box}_1),\ldots, c({\rm box}_n))$ is 
$(-7,-6,-5,-2,0,1,1,2,2,3,3,4,5,6)$.
$$
\matrix{
\beginpicture
\setcoordinatesystem units <0.5cm,0.5cm>         
\setplotarea x from 0 to 4, y from 0 to 3    
\linethickness=0.5pt                          
\put{-7} at 0.5 0.5
\put{-6} at 0.5 1.5
\put{-5} at 1.5 1.5
\put{-2}  at 3.5 2.5
\put{0} at 4.5 3.5
\put{1} at 4.5 4.5
\put{1}  at 5.5 3.5
\put{2}  at 5.5 4.5
\put{3}  at 5.5 5.5
\put{2}  at 6.5 3.5
\put{3}  at 6.5 4.5
\put{4}  at 6.5 5.5
\put{5}  at 7.5 5.5
\put{6}  at 8.5 5.5
\putrule from 5 6 to 9 6          %
\putrule from 4 5 to 9 5          
\putrule from 4 4 to 7 4          %
\putrule from 3 3 to 7 3          %
\putrule from 3 2 to 4 2          %
\putrule from 0 2 to 2 2          %
\putrule from 0 1 to 2 1          %
\putrule from 0 0 to 1 0          %
\putrule from 0 0 to 0 2        %
\putrule from 1 0 to 1 2        %
\putrule from 2 1 to 2 2        %
\putrule from 3 2 to 3 3        
\putrule from 4 2 to 4 5        %
\putrule from 5 3 to 5 6        %
\putrule from 6 3 to 6 6        %
\putrule from 7 3 to 7 6        %
\putrule from 8 5 to 8 6        %
\putrule from 9 5 to 9 6        %
\endpicture
\cr
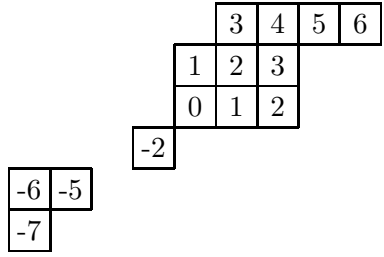
\hbox{Contents of the boxes of $(c,\lambda/\mu)$} \cr
}
$$
The following picture shows the contents of the boxes in the placed skew
shape $(c',\lambda/\mu)$ such that $(c({\rm
box}_1),\ldots, c({\rm
box}_n))=(-7,-6,-5,-3/2,1/2,3/2,3/2,5/2,5/2,7/2,7/2,9/2,11/2,13/2)$.
$$
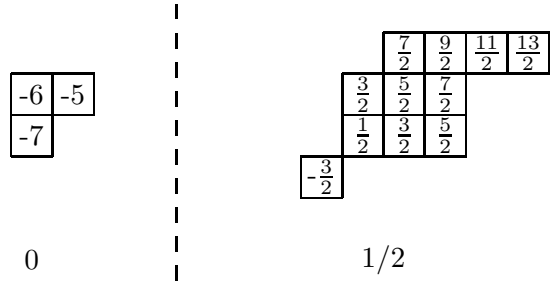
\beginpicture
\setcoordinatesystem units <0.5cm,0.5cm>         
\setplotarea x from -4 to 4, y from 0 to 7    
\linethickness=0.5pt                          
\put{-7} at -3.5 3.5
\put{-6} at -3.5 4.5
\put{-5} at -2.5 4.5
\put{-${3\over2}$}  at 3.5 2.5
\put{$1\over2$} at 4.5 3.5
\put{$3\over2$} at 4.5 4.5
\put{$3\over2$}  at 5.5 3.5
\put{$5\over 2$}  at 5.5 4.5
\put{$7\over 2$}  at 5.5 5.5
\put{$5\over 2$}  at 6.5 3.5
\put{$7\over2$}  at 6.5 4.5
\put{$9\over2$}  at 6.5 5.5
\put{$11\over2$}  at 7.5 5.5
\put{$13\over2$}  at 8.5 5.5
\putrule from -4 5 to -2 5          %
\putrule from -4 4 to -2 4          
\putrule from -4 3 to -3 3          %
\putrule from -4 3 to -4 5        %
\putrule from -3 3 to -3 5        
\putrule from -2 4 to -2 5        %
\putrule from 5 6 to 9 6          %
\putrule from 4 5 to 9 5          
\putrule from 4 4 to 7 4          
\putrule from 3 3 to 7 3          %
\putrule from 3 2 to 4 2          %
\putrule from 3 2 to 3 3        
\putrule from 4 2 to 4 5        
\putrule from 5 3 to 5 6        %
\putrule from 6 3 to 6 6        %
\putrule from 7 3 to 7 6        %
\putrule from 8 5 to 8 6        %
\putrule from 9 5 to 9 6        %
\setdashes
\putrule from 0 0 to 0 7           
\put{0}   at -3.5 0.5      %
\put{1/2} at 5 0.5       
\endpicture
$$
This ``book'' has two pages, with page numbers $0$ and $1/2$.
\endpf

\lemma  Let $(c,\lambda/\mu)$ be a placed skew shape with $n$ boxes and let $L$
be a  standard tableau of shape $\lambda/\mu$.  Let $L(i)$ denote the
box containing $i$ in $L$.
The {\it content sequence}
$$(c(L(1)),\ldots,c(L(n))$$
uniquely determines the shape $(c,\lambda/\mu)$ and the standard
tableau $L$.
\pf
Proceed by induction on the number of boxes of $L$.
If $L$ has only one box then the content sequence $(c(L(1))$
determines the placement of that box.
Assume that $L$ has $n$ boxes.  Let $L'$ be the standard tableau
determined by removing the box containing $n$ from $L$.
Then $L'$ is also of skew shape and the content sequence of
$L'$ is $(c(L(1)),\ldots,c(L(n-1)))$.  By the induction hypothesis
we can reconstruct $L'$ from its content sequence.  Then $c(L(n))$ determines
the diagonal which must contain box $n$ in $L$.  So $L'$ and $c(L(n))$
determine $L$ uniquely.
\endpf

\bigskip

\section 3. Weights and weight spaces

\bigskip

A finite dimensional $\tilde H_n$-module is {\it calibrated}
if it has a basis of simultaneous eigenvectors for the $x_i$,
$1\le i\le n$.    In other words, $M$ is calibrated if it has a basis
$\{v_t\}$ such that  for all $v_t$ in the basis and all $1\le i\le n$,
$$x_iv_t = t_i v_t,
\qquad \hbox{for some $t_i\in \CC^*$.}$$

\subsection{Weights.}

Let $X$ be the abelian group generated by the elements $x_1,\ldots, x_n\in
\tilde H_n$ and let 
$$T=\{ \hbox{group homomorphisms $X\to \CC^*$}\}.$$
The {\it torus} $T$ can be identified with $(\CC^*)^n$ by
identifying the element $t=(t_1,\ldots, t_n)\in (\CC^*)^n$ with the homomorphism
given by $t(x_i) = t_i$, $1\le i\le n$. The symmetric group $S_n$ acts on
$\ZZ^n$ by permuting coordinates and this action induces an action
of $S_n$ on $T$ given by
$$(wt)(x^\gamma)= x^{w^{-1}\gamma},
\qquad\hbox{for $w\in S_n$, $\gamma\in \ZZ^n$,}$$
with notation as in (1.2).

\subsection{Weight spaces.}

Let $M$ be a finite dimensional $\tilde H_n$-module.
For each $t=(t_1,\ldots, t_n)\in T$ the {\it $t$-weight space of $M$} 
and the {\it generalized $t$-weight space} are the subspaces
$$\eqalign{
M_t &= \{ m\in M \ |\ x_i m = t_i m
\hbox{\ for all $1\le i\le n$}\}
\qquad\hbox{and}  \cr
\cr
M_t^{\rm gen} &= \{ m\in M\ |\
\hbox{for each $1\le i\le n$,
$(x_i-t_i)^k m=0$ for some $k\in \ZZ_{>0}$}
\},  \cr
}
$$
respectively.
From the definitions, $M_t\subseteq M_t^{\rm gen}$ and $M$ is calibrated if
and only if
$M_t^{\rm gen}=M_t$ for all $t\in T$.   If $M_t^{\rm gen}\ne 0$ then $M_t\ne 0$. 
In general $M\neq\bigoplus_{t\in T} M_t$, but we do have
$$
M=\bigoplus_{t\in T} M_t^{\rm gen}.
$$
This is a decomposition of $M$ into Jordan blocks
for the action of $\CC[X]=\CC[x_1^{\pm 1},\ldots, x_n^{\pm 1}]$.  
The set of {\it weights of $M$} is the set
$$\supp(M) = \{ t\in T \ |\ M_t^{\rm gen}\ne 0\}.\formula$$
An element of $M_t$ is called a {\it weight vector} of {\it weight} $t$.

\subsection{The $\tau$ operators.}

The maps $\tau_i\colon M_t^{\rm gen}\to M_{s_it}^{\rm gen}$ 
defined below are local operators on $M$ in the sense that they
act on each generalized weight space $M_t^{\rm gen}$ of $M$ separately.  The
operator $\tau_i$ is only defined on the generalized weight spaces $M_t^{\rm
gen}$ such that $t_i\ne t_{i+1}$.  

\prop 
Let $t=(t_1,\ldots,t_n)\in T$ be such that $t_i\ne t_{i+1}$ and let $M$
be a finite dimensional $\tilde H_n$-module.  Define
$$
\matrix{
\tau_i \colon 
&M_t^{\rm gen} &\longrightarrow &M_{s_it}^{\rm gen} \cr
\cr
& m & \longmapsto &
\displaystyle{
\left(T_i - {(q-q^{-1})x_{i+1}\over x_{i+1}-x_i }\right) m.} \cr
}
$$
\item{(a)}
The map $\tau_i\colon M_t^{\rm gen} \longrightarrow M_{s_it}^{\rm gen}$
is well defined.
\smallskip
\item{(b)}  As operators on $M_t^{\rm gen}$
$$\eqalign{
x_i\tau_i &= \tau_i x_{i+1}, \qquad
x_{i+1}\tau_i = \tau_i x_i\;,  \qquad\hbox{and}\qquad x_j\tau_i = \tau_ix_j,
\quad\hbox{if $j\ne i,i+1$}, \cr
\tau_i\tau_i&=
{(qx_{i+1}-q^{-1}x_i)(qx_i-q^{-1}x_{i+1})\over 
(x_{i+1}-x_i)(x_i-x_{i+1}) },
\qquad\qquad\quad
\hbox{if\ \ $1\le i\le n-1$,} \cr
\tau_i\tau_j &= \tau_j\tau_i, 
\qquad\qquad\qquad\qquad\qquad\qquad\qquad\qquad\qquad
\hbox{if\ \  $|i-j|>1$,}
\cr
\tau_i\tau_{i+1}\tau_i &= \tau_{i+1}\tau_i\tau_{i+1},
\qquad\quad\qquad\qquad\qquad\qquad\qquad\qquad\   1\le i\le n-1,  \cr
}$$
whenever both sides are well defined.
\pf
(a)  Note that $(q-q^{-1})x_{i+1}/(x_{i+1}-x_i)$ is not a well defined
element of $\tilde H_n$ or $\CC[x_1^{\pm1},\ldots,x_n^{\pm1}]$
since it is a power series and not a Laurent polynomial.  Because of this we
will be careful to view $(q-q^{-1})x_{i+1}/(x_{i+1}-x_i)$ only as an
operator on $M_t^{\rm gen}$.  Let us describe this operator more precisely.

The element $x_ix_{i+1}^{-1}$ acts on $M_t^{\rm gen}$ by $t_it_{i+1}^{-1}$
times a unipotent transformation.  As an operator on $M_t^{\rm gen}$,
$(1-x_ix_{i+1}^{-1}) = x_{i+1}/(x_{i+1}-x_i)$ is invertible since it has 
determinant $(1-t_it_{i+1}^{-1})^d$ where $d=\dim(M_t^{\rm gen})$.  Since this
determinant is nonzero $(q-q^{-1})x_{i+1}/(x_{i+1}-x_i)
=(q-q^{-1})(1-x_ix_{i+1}^{-1})^{-1}$ is a well defined operator on $M_t^{\rm
gen}$.  Thus the definition of $\tau_i$ makes sense.

The operator identities 
$x_i\tau_i = \tau_i x_{i+1}$
$x_{i+1}\tau_i = \tau_i x_i$,  and $x_j\tau_j = \tau_ix_j$,
if $i\ne i,i+1$, now follow easily from the definition of the $\tau_i$ and the
identities in (1.4).  These identities imply that 
$\tau_i$ maps $M_t^{\rm gen}$ into $M_{s_it}^{\rm gen}$.  

All of the operator identities in part (b) are proved by straightforward
calculations of the same flavour as the calculation of $\tau_i\tau_i$ given
below.  We shall not give the details for the other cases.
The only one which is really tedious is the calculation for the proof of
$\tau_{i+1}\tau_i\tau_{i+1} =\tau_i\tau_{i+1}\tau_i$.  For a
more pleasant (but less elementary) proof of this identity see Proposition 2.7
in [Ra3].

Since $t_i\ne t_{i+1}$ both $\tau_i\colon M_t^{\rm gen}\to M_{s_it}^{\rm gen}$
and $\tau_i\colon M_{s_it}^{\rm gen}\to M_t^{\rm gen}$ are well defined.
Let $m\in M_t^{\rm gen}$.  Then
$$\eqalign{
\tau_i\tau_im
&=
\left(T_i - {(q-q^{-1})x_{i+1}\over x_{i+1}-x_i }\right)
\left(T_i - {(q-q^{-1})x_{i+1}\over x_{i+1}-x_i }\right) m \cr
&=
\left(T_i^2 
-T_i {(q-q^{-1})x_{i+1}\over x_{i+1}-x_i }
- {(q-q^{-1})x_{i+1}\over x_{i+1}-x_i }T_i
+{(q-q^{-1})^2x_{i+1}^2\over (x_{i+1}-x_i)^2 }
\right) m \cr
&=
\left((q-q^{-1})T_i +1 - T_i {(q-q^{-1})x_{i+1}\over x_{i+1}-x_i }
-T_i {(q-q^{-1})x_i\over x_i-x_{i+1} } \right. \cr
&\qquad\qquad
\left.
-(q-q^{-1})^2 {x_{i+1}\over x_{i+1}-x_i}
\right.\left( {x_{i+1}\over x_{i+1}-x_i} - {x_i\over x_i-x_{i+1}}
\right)\left.
+{(q-q^{-1})^2x_{i+1}^2\over (x_{i+1}-x_i)^2 }
\right) m \cr
&=
\left((q-q^{-1})T_i +1 - (q-q^{-1})T_i
+(q-q^{-1})^2{x_ix_{i+1}\over (x_{i+1}-x_i)(x_i-x_{i+1})}
\right) m \cr
&=
 {q^2x_ix_{i+1}-2x_ix_{i+1}+q^{-2}x_ix_{i+1}-x_i^2+2x_ix_{i+1}-x_{i+1}^2\over 
(x_{i+1}-x_i)(x_i-x_{i+1}) }
 m \cr
&=
{ (qx_{i+1}-q^{-1}x_i)(qx_i-q^{-1}x_{i+1})\over 
(x_{i+1}-x_i)(x_i-x_{i+1}) } 
 m. \qquad\qquad\qquad\qquad\hbox{\qed}\cr
}
$$

Let $w\in S_n$.  Let $w=s_{i_1}\cdots s_{i_p}$ be a reduced word
for $w$ and define
$$\tau_w = \tau_{i_1}\cdots \tau_{i_p}.\formula$$
Since the $\tau$-operators satisfy the braid relations the operator $\tau_w$
is independent of the choice of the reduced word for $w$.
The operator $\tau_w$ is a well defined operator on $M_t^{\rm gen}$
if $t=(t_1,\ldots,t_n)$ is such that $t_i\ne t_j$ for all
pairs $i<j$ such that $w(i)>w(j)$.
One may use the relations in (1.5) to rewrite $\tau_w$ in the form
$$\tau_w = \sum_{u\le w} T_w a_{uw}(x_1,\ldots,x_n)$$ 
where $a_{uw}(x_1,\ldots, x_n)$ are rational functions
in the variables $x_1,\ldots, x_n$. (The functions $a_{uw}(x_1,\ldots,x_n)$
are analogues of the Harish-Chandra {\it $c$-function}, see
[Mac2, 4.1] and [Op, Theorem 5.3].) 
If $t=(t_1,\ldots,t_n)$ is such that $t_i\ne t_j$ for all
pairs $i<j$ such that $w(i)>w(j)$ then the expression
$$\tau_w\big\vert_t = \sum_{u\le w} T_w a_{uw}(t_1,\ldots,t_n)\formula$$
is a well defined element of the Iwahori-Hecke algebra $H_n$.
If $w=uv$ with $\ell(w)=\ell(u)+\ell(v)$ then
$$\tau_w \big\vert_t = \tau_u\big\vert_{vt} \tau_v\big\vert_t.\formula$$
The following result will be crucial to the proof of Theorem 5.5.
This result is due to D. Barbasch and P. Diaconis [D] (in the $q=1$ case). 
The proof given below is a $q$-version of a proof for the $q=1$ case given by S.
Fomin [Fo].

\prop  Let $w_0$ be the longest element of $S_n$,
$$w_0 = \pmatrix{1 &2 &\cdots &n-1 &n\cr
n &n-1 &\cdots &2 &1}.$$ 
Let $a\in \CC^*$ and fix $t=(a,aq^2,aq^4,\ldots, aq^{2(n-1)})$. Then
$$\tau_{w_0}\big\vert_t = 
\sum_{w\in S_n} T_w (-q)^{\ell(w_0)-\ell(w)},$$
where $\ell(w_0)={n\choose 2}$.
\pf
Let $1\le k\le n$. Then there is a $v\in S_n$ such that
$w_0 = vs_k$ and $\ell(w_0) = \ell(v)+1$.  
So
$$\tau_{w_0} = \tau_v\tau_k
=\tau_v\left(T_k - {(q-q^{-1})x_{k+1}\over x_{k+1}-x_k}\right)$$
and
$$
\tau_{w_0}\big\vert_t
= \tau_v\big\vert_{s_kt}
\left(T_k - {(q-q^{-1})t_{k+1}\over t_{k+1}-t_k}\right)
= \tau_v\big\vert_{s_kt}
\left(T_k - {(q-q^{-1})q^2t_k\over (q^2-1)t_k}\right)
= \tau_v\big\vert_{s_kt}
(T_k - q).
$$
Right multiplying by $T_k+q^{-1}$ and using the relation (1.3) gives
$$
\tau_{w_0}\big\vert_t (T_k +q^{-1}) = \tau_v\big\vert_{s_kt} 
(T_k - q)(T_k+q^{-1})
= 0.$$
The element 
$h = \sum_{w\in S_n} T_w (-q)^{\ell(w_0)-\ell(w)}$
is a multiple of the minimal central idempotent in $H_n$ 
corresponding to the representation $\phi$ given by
$\phi(T_k) = -q^{-1}$, for all $1\le k\le n$.  Up to multiplication by
constants, it is the unique element in $H_n$ such that $h(T_k+q^{-1})=0$
for all $1\le k\le n$.  The lemma follows by noting that the
coefficients of $T_{w_0}$ in $h$ and $\tau_{w_0}\big\vert_t$ are both $1$.
\endpf

The action of the $\tau$-operators on weight vectors will be particularly
important to the proofs of the results in later sections.
Let us record the following facts.
\smallskip\noindent
Let $M$
be an $\tilde H_n$-module and let $m_t$ be a weight vector in $M$
of weight $t=(t_1,\ldots,t_n)$.
\smallskip
\item{(\global\advance\resultno by 1
\the\sectno.\the\resultno a)}  If $t_i\ne t_{i+1}$ then
$$\tau_i m_t = \left(T_i - {(q-q^{-1})x_{i+1}\over x_i-x_{i+1}}\right)m_t
=\left(T_i - {(q-q^{-1})t_{i+1}\over t_i-t_{i+1}}\right)m_t$$
is a weight vector of weight $s_it$.  
\smallskip
\item{(\the\sectno.\the\resultno b)}   By the second set of identities in 
Proposition 3.2 (b),
$\tau_i\tau_i m_t = (qt_{i+1}-q^{-1}t_i)(qt_i-q^{-1}t_{i+1})
(t_{i+1}-t_i)^{-1}(t_i-t_{i+1})^{-1}m_t$.  Thus
$$\hbox{If $t_i\ne t_{i+1}$ and $t_i\ne q^{\pm
2}t_{i+1}$ then
$\tau_i m_t\ne 0$.}$$

\vfill\eject

\section 4.  Classification and construction of calibrated representations

\bigskip

The following theorem classifies and constructs all irreducible calibrated
representations of the affine Hecke algebra $\tilde H_n$.  It shows that the
theory of standard Young tableaux plays an intrinsic role in the combinatorics
of the representations of the affine Hecke algebra.  The construction given in
Theorem 4.1 is a direct generalization of A. Young's classical ``seminormal
construction'' of the irreducible representations of the symmetric group [Y].
Young's construction was generalized to Iwahori-Hecke algebras of type A by
Hoefsmit [Ho] and Wenzl [Wz] independently, to Iwahori-Hecke algebras of types
B and D by Hoefsmit [Ho] and to cyclotomic Hecke algebras by Ariki and Koike
[AK].  It can be shown that all of these previous generalizations are special
cases of the construction for affine Hecke algebras given here.  Recently, this
construction has been generalized even further [Ra3], to affine Hecke algebras of
arbitrary Lie type.  Some parts of Theorem 4.1 are due, originally, to
I. Cherednik, and are stated in [Ch, \S 3].

Garsia and Wachs [GW] showed 
that the theory of standard Young tableaux  and Young's constructions play an
important role in the combinatorics of the skew representations of the symmetric
group.  At that time it was not known that these
representations are actually {\it irreducible} as representations of the affine
Hecke algebra!!

\thm  Let $(c,\lambda/\mu)$ be a placed skew shape with $n$ boxes.  Define an
action of $\tilde H_n$ on the vector space
$$\tilde H^{(c,\lambda/\mu)} =
\CC{\rm -span}\{ v_L\ |\ \hbox{$L$ is a standard
tableau of shape $\lambda/\mu$}\}$$ 
by the formulas
$$\eqalign{
x_i v_L &= q^{2c(L(i))} v_L, \cr
T_i v_L &= (T_i)_{LL} v_L + (q^{-1}+(T_i)_{LL}) v_{s_iL}, \cr
}
$$
where $s_iL$ is the same as $L$ except that the entries $i$ and $i+1$ are
interchanged,
$$(T_i)_{LL} = { q-q^{-1} \over 1-q^{2(c(L(i))-c(L(i+1)))} },
\quad\quad
\hbox{$v_{s_iL}=0$ if $s_iL$ is not a standard tableau,}$$
and $L(i)$ denotes the box of $L$ containing the entry $i$.
\smallskip
\item{(a)}  $\tilde H^{(c,\lambda/\mu)}$ is a calibrated irreducible
$\tilde H_n$-module.
\smallskip
\item{(b)}  The modules $\tilde H^{(c,\lambda/\mu)}$ are
non-isomorphic.
\smallskip
\item{(c)}  Every irreducible calibrated $\tilde H_n$-module is
isomorphic to $\tilde H^{(c,\lambda/\mu)}$ for some placed skew shape
$(c,\lambda/\mu)$.
\endthm  
\smallskip\noindent
{\smallcaps step 1.} The given formulas for the action of $\tilde
H^{(c,\lambda/\mu)}$ define an $\tilde H_n$-module.

\pf
If $L$ is a standard tableau then the entries $i$ and $i+1$ cannot
appear in the same diagonal in $L$.  Thus, for all standard tableaux $L$,
$c(L(i))\ne c(L(i+1))$ and for this reason the constant
$(T_i)_{LL}$ is always well defined. 

Let $L$ be a standard tableau of shape $\lambda/\mu$.  Then
$(T_i)_{LL}+(T_i)_{s_iL,s_iL}=q-q^{-1}$ and so
$$\eqalign{
T_i^2v_L 
&= ((T_i)_{LL}^2 
+ (q^{-1}+(T_i)_{LL})(q^{-1}+(T_i)_{s_iL,s_iL}))v_L \cr
&\qquad\qquad\quad
+ (q^{-1}+(T_i)_{LL})((T_i)_{LL}+(T_i)_{s_iL,s_iL}) v_{s_iL} 
\cr
&= (T_i)_{LL}( (T_i)_{LL}+(T_i)_{s_iL,s_iL} ) v_L
+q^{-1}(q^{-1} + (T_i)_{LL}+(T_i)_{s_iL,s_iL}) )v_L \cr
&\qquad\qquad\quad
+ (q^{-1}+(T_i)_{LL}) (q-q^{-1}) v_{s_iL}\cr
&= (T_i)_{LL}(q-q^{-1}) v_L +(q^{-1}+(T_i)_{LL})(q-q^{-1})v_{s_iL}
+q^{-1}(q^{-1}+q-q^{-1}) v_L \cr
&= ((q-q^{-1})T_i +1)v_L.\cr
}
$$
The calculations to check the identities (1.1a), (1.1f) and (1.5) are routine.
Checking the identity $T_iT_{i+1}T_i=T_{i+1}T_iT_{i+1}$ is more involved. One
can proceed as follows.  According to the formulas for the action, the operators
$T_i$ and $T_{i+1}$ preserve the subspace $S$ spanned by the vectors $v_Q$
indexed by the standard tableaux $Q$ in the set
$\{L, s_iL, s_{i+1}L, s_is_{i+1}L, s_{i+1}s_iL,
s_is_{i+1}s_iL\}$.  Depending on the relative positions of the boxes containing
$i, i+1, i+2$ in $L$, this space is either $1,2,3$ or $6$ dimensional. 
Representative cases are when these boxes are positioned in the following ways. 
$$\matrix
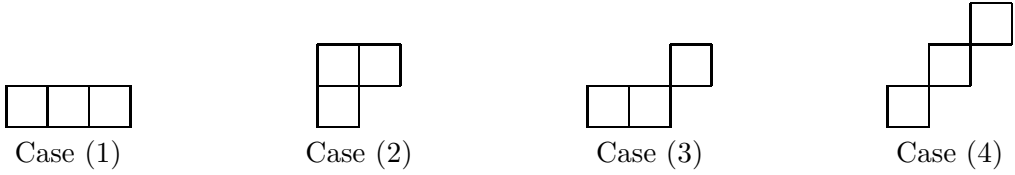
{
\beginpicture
\setcoordinatesystem units <0.5cm,0.5cm>         
\setplotarea x from 0 to 3, y from 0 to 3    
\linethickness=0.5pt                          
\putrule from 0 1 to 3 1          %
\putrule from 0 0 to 3 0          
\putrule from 0 0 to 0 1        %
\putrule from 1 0 to 1 1        %
\putrule from 2 0 to 2 1        
\putrule from 3 0 to 3 1        %
\endpicture
&\qquad\qquad
&\beginpicture
\setcoordinatesystem units <0.5cm,0.5cm>         
\setplotarea x from 0 to 2, y from 0 to 3    
\linethickness=0.5pt                          
\putrule from 0 2 to 2 2          %
\putrule from 0 1 to 2 1          %
\putrule from 0 0 to 1 0          
\putrule from 0 0 to 0 2        %
\putrule from 1 0 to 1 2        %
\putrule from 2 1 to 2 2        
\endpicture
&\qquad\qquad
&\beginpicture
\setcoordinatesystem units <0.5cm,0.5cm>         
\setplotarea x from 0 to 3, y from 0 to 3    
\linethickness=0.5pt                          
\putrule from 0 0 to 2 0          %
\putrule from 0 1 to 3 1          %
\putrule from 2 2 to 3 2          
\putrule from 0 0 to 0 1        %
\putrule from 1 0 to 1 1        %
\putrule from 2 0 to 2 2        
\putrule from 3 1 to 3 2        %
\endpicture
&\qquad\qquad
&\beginpicture
\setcoordinatesystem units <0.5cm,0.5cm>         
\setplotarea x from 0 to 3, y from 0 to 3    
\linethickness=0.5pt                          
\putrule from 0 0 to 1 0          %
\putrule from 0 1 to 2 1          %
\putrule from 1 2 to 3 2          
\putrule from 2 3 to 3 3          
\putrule from 0 0 to 0 1        %
\putrule from 1 0 to 1 2        %
\putrule from 2 1 to 2 3        
\putrule from 3 2 to 3 3        %
\endpicture  \cr
\hbox{Case (1)}
&&\hbox{Case (2)}
&&\hbox{Case (3)}
&&\hbox{Case (4)} \cr
}
$$
In Case (1) the space $S$ is one dimensional and spanned by the vector
$v_Q$ corresponding to the standard tableau
$$\beginpicture
\setcoordinatesystem units <0.5cm,0.5cm>         
\setplotarea x from 0 to 3, y from 0 to 1.5   
\linethickness=0.5pt                          
\put{$a$} at .5 .5 
\put{$b$} at 1.5 .5
\put{$c$} at 2.5 .5
\putrule from 0 1 to 3 1          %
\putrule from 0 0 to 3 0          
\putrule from 0 0 to 0 1        %
\putrule from 1 0 to 1 1        %
\putrule from 2 0 to 2 1        
\putrule from 3 0 to 3 1        %
\endpicture
$$
where $a=i$, $b=i+1$, and $c=i+2$.  The action
of $T_i$ and $T_{i+1}$ on $S$ is given by the matrices
$$\phi_S(T_i)=\pmatrix{q},\qquad\hbox{and}\qquad
\phi_S(T_{i+1})=\pmatrix{q},$$
respectively.  In case (2) the space $S$ is two dimensional and spanned by the
vectors $v_Q$ corresponding to the standard tableaux
$$\beginpicture
\setcoordinatesystem units <0.5cm,0.5cm>         
\setplotarea x from 0 to 2, y from 0 to 2    
\linethickness=0.5pt                          
\put{$a$} at .5 1.5
\put{$b$} at 1.5 1.5
\put{$c$} at .5 .5
\putrule from 0 2 to 2 2          %
\putrule from 0 1 to 2 1          %
\putrule from 0 0 to 1 0          
\putrule from 0 0 to 0 2        %
\putrule from 1 0 to 1 2        %
\putrule from 2 1 to 2 2        
\endpicture
\qquad\qquad\qquad\qquad
\beginpicture
\setcoordinatesystem units <0.5cm,0.5cm>         
\setplotarea x from 0 to 2, y from 0 to 2    
\linethickness=0.5pt                          
\put{$a$} at .5 1.5
\put{$c$} at 1.5 1.5
\put{$b$} at .5 .5
\putrule from 0 2 to 2 2          %
\putrule from 0 1 to 2 1          %
\putrule from 0 0 to 1 0          
\putrule from 0 0 to 0 2        %
\putrule from 1 0 to 1 2        %
\putrule from 2 1 to 2 2        
\endpicture
$$
where $a=i$, $b=i+1$, and $c=i+2$.  The action of $T_i$ and
$T_{i+1}$ on
$S$ is given by the matrices
$$\phi_S(T_{i})=\pmatrix{q &0 \cr 0 &-q^{-1} \cr}
\qquad\hbox{and}\qquad
\phi_S(T_{i+1})=\pmatrix{
\displaystyle{ {q-q^{-1}\over 1-q^4} } 
&\displaystyle{ {q-q^{-5}\over 1-q^{-4}} }     \cr 
\displaystyle{ {q-q^3\over 1-q^4} }    
&\displaystyle{ {q-q^{-1}\over 1-q^{-4}} } \cr
}
$$
In case (3) the space $S$ is three dimensional and 
spanned by the vectors $v_Q$ corresponding to the standard tableaux
$$
\beginpicture
\setcoordinatesystem units <0.5cm,0.5cm>         
\setplotarea x from 0 to 3, y from 0 to 2   
\linethickness=0.5pt                          
\put{$a$} at .5 .5
\put{$b$} at 1.5 .5
\put{$c$} at 2.5 1.5
\putrule from 0 0 to 2 0          %
\putrule from 0 1 to 3 1          %
\putrule from 2 2 to 3 2          
\putrule from 0 0 to 0 1        %
\putrule from 1 0 to 1 1        %
\putrule from 2 0 to 2 2        
\putrule from 3 1 to 3 2        %
\endpicture
\qquad\qquad
\beginpicture
\setcoordinatesystem units <0.5cm,0.5cm>         
\setplotarea x from 0 to 3, y from 0 to 2    
\linethickness=0.5pt                          
\put{$a$} at .5 .5
\put{$c$} at 1.5 .5
\put{$b$} at 2.5 1.5
\putrule from 0 0 to 2 0          %
\putrule from 0 1 to 3 1          %
\putrule from 2 2 to 3 2          
\putrule from 0 0 to 0 1        %
\putrule from 1 0 to 1 1        %
\putrule from 2 0 to 2 2        
\putrule from 3 1 to 3 2        %
\endpicture
\qquad\qquad
\beginpicture
\setcoordinatesystem units <0.5cm,0.5cm>         
\setplotarea x from 0 to 3, y from 0 to 2    
\linethickness=0.5pt                          
\put{$b$} at .5 .5
\put{$c$} at 1.5 .5
\put{$a$} at 2.5 1.5
\putrule from 0 0 to 2 0          %
\putrule from 0 1 to 3 1          %
\putrule from 2 2 to 3 2          
\putrule from 0 0 to 0 1        %
\putrule from 1 0 to 1 1        %
\putrule from 2 0 to 2 2        
\putrule from 3 1 to 3 2        %
\endpicture
$$
where $a=i$, $b=i+1$, and $c=i+2$.  The action of
$T_i$ and $T_{i+1}$ on $S$ is given by the matrices
$$\eqalign{
\phi_S(T_{i})&=\pmatrix{q &0 &0\cr 
0 
&\displaystyle{ {q-q^{-1}\over 1-q^{2(c_1-c_3)}} } 
&\displaystyle{ {q-q^{2(c_3-c_1)-1}\over 1-q^{2(c_3-c_1)}} } \cr 
0 
&\displaystyle{ {q-q^{2(c_1-c_3)-1}\over 1-q^{2(c_1-c_3)}} }  
&\displaystyle{ {q-q^{-1}\over 1-q^{2(c_3-c_1)}} } \cr}
\qquad\hbox{and} \cr \cr
\phi_S(T_{i+1})&=\pmatrix{
\displaystyle{ {q-q^{-1}\over 1-q^{2(c_2-c_3)}} } 
&\displaystyle{ {q-q^{2(c_3-c_2)-1}\over 1-q^{2(c_3-c_2)}} }
&0 \cr
\displaystyle{ {q-q^{2(c_2-c_3)-1}\over 1-q^{2(c_2-c_3)}} }
&\displaystyle{ {q-q^{-1}\over 1-q^{2(c_3-c_2)}} } 
&0 \cr
0  &0  &q \cr
} \cr
}
$$
respectively, where $c_1 = c(L(i))$, $c_2=c(L(i+1))$ and
$c_3=c(L(i+2))$.  In case (4) the space
$S$ is six dimensional and spanned by the
vectors $v_Q$ corresponding to the standard tableaux
$$
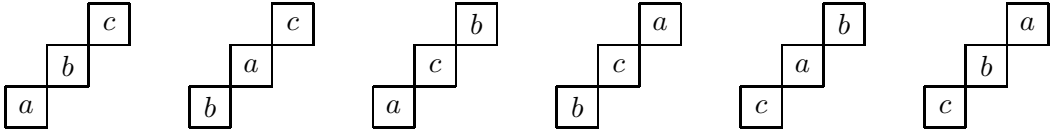
\beginpicture
\setcoordinatesystem units <0.5cm,0.5cm>         
\setplotarea x from 0 to 3, y from 0 to 3    
\linethickness=0.5pt                          
\put{$a$} at .5 .5
\put{$b$} at 1.5 1.5
\put{$c$} at 2.5 2.5
\putrule from 0 0 to 1 0          %
\putrule from 0 1 to 2 1          %
\putrule from 1 2 to 3 2          
\putrule from 2 3 to 3 3          
\putrule from 0 0 to 0 1        %
\putrule from 1 0 to 1 2        %
\putrule from 2 1 to 2 3        
\putrule from 3 2 to 3 3        %
\endpicture  
\qquad
\beginpicture
\setcoordinatesystem units <0.5cm,0.5cm>         
\setplotarea x from 0 to 3, y from 0 to 3    
\linethickness=0.5pt                          
\put{$b$} at .5 .5
\put{$a$} at 1.5 1.5
\put{$c$} at 2.5 2.5
\putrule from 0 0 to 1 0          %
\putrule from 0 1 to 2 1          %
\putrule from 1 2 to 3 2          
\putrule from 2 3 to 3 3          
\putrule from 0 0 to 0 1        %
\putrule from 1 0 to 1 2        %
\putrule from 2 1 to 2 3        
\putrule from 3 2 to 3 3        %
\endpicture  
\qquad
\beginpicture
\setcoordinatesystem units <0.5cm,0.5cm>         
\setplotarea x from 0 to 3, y from 0 to 3    
\linethickness=0.5pt                          
\put{$a$} at .5 .5
\put{$c$} at 1.5 1.5
\put{$b$} at 2.5 2.5
\putrule from 0 0 to 1 0          %
\putrule from 0 1 to 2 1          %
\putrule from 1 2 to 3 2          
\putrule from 2 3 to 3 3          
\putrule from 0 0 to 0 1        %
\putrule from 1 0 to 1 2        %
\putrule from 2 1 to 2 3        
\putrule from 3 2 to 3 3        %
\endpicture  
\qquad
\beginpicture
\setcoordinatesystem units <0.5cm,0.5cm>         
\setplotarea x from 0 to 3, y from 0 to 3    
\linethickness=0.5pt                          
\put{$b$} at .5 .5
\put{$c$} at 1.5 1.5
\put{$a$} at 2.5 2.5
\putrule from 0 0 to 1 0          %
\putrule from 0 1 to 2 1          %
\putrule from 1 2 to 3 2          
\putrule from 2 3 to 3 3          
\putrule from 0 0 to 0 1        %
\putrule from 1 0 to 1 2        %
\putrule from 2 1 to 2 3        
\putrule from 3 2 to 3 3        %
\endpicture  
\qquad
\beginpicture
\setcoordinatesystem units <0.5cm,0.5cm>         
\setplotarea x from 0 to 3, y from 0 to 3    
\linethickness=0.5pt                          
\put{$c$} at .5 .5
\put{$a$} at 1.5 1.5
\put{$b$} at 2.5 2.5
\putrule from 0 0 to 1 0          %
\putrule from 0 1 to 2 1          %
\putrule from 1 2 to 3 2          
\putrule from 2 3 to 3 3          
\putrule from 0 0 to 0 1        %
\putrule from 1 0 to 1 2        %
\putrule from 2 1 to 2 3        
\putrule from 3 2 to 3 3        %
\endpicture  
\qquad
\beginpicture
\setcoordinatesystem units <0.5cm,0.5cm>         
\setplotarea x from 0 to 3, y from 0 to 3    
\linethickness=0.5pt                          
\put{$c$} at .5 .5
\put{$b$} at 1.5 1.5
\put{$a$} at 2.5 2.5
\putrule from 0 0 to 1 0          %
\putrule from 0 1 to 2 1          %
\putrule from 1 2 to 3 2          
\putrule from 2 3 to 3 3          
\putrule from 0 0 to 0 1        %
\putrule from 1 0 to 1 2        %
\putrule from 2 1 to 2 3        
\putrule from 3 2 to 3 3        %
\endpicture  
$$
where $a=i$, $b=i+1$, and $c=i+2$.  The action of $T_i$ and $T_{i+1}$ on $S$ is
given by the matrices
$$
\phi_S(T_{i})=\pmatrix{
\displaystyle{ {q-q^{-1}\over 1-q^{2d_{12}}} } 
&\displaystyle{ {q-q^{2d_{21}-1}\over 1-q^{2d_{21}}} }
&0 &0 &0 &0\cr
\displaystyle{ {q-q^{2d_{12}-1}\over 1-q^{2d_{12}}} }
&\displaystyle{ {q-q^{-1}\over 1-q^{2d_{21}}} } 
&0 &0 &0 &0 \cr
0 &0 &\displaystyle{ {q-q^{-1}\over 1-q^{2d_{13}}} } 
&\displaystyle{ {q-q^{2d_{31}-1}\over 1-q^{2d_{31}}} }
&0 &0 \cr
0 &0 &\displaystyle{ {q-q^{2d_{13}-1}\over 1-q^{2d_{13}}} }
&\displaystyle{ {q-q^{-1}\over 1-q^{2d_{31}}} } 
&0 &0 \cr
0 &0 &0 &0 &\displaystyle{ {q-q^{-1}\over 1-q^{2d_{23}}} } 
&\displaystyle{ {q-q^{2d_{32}-1}\over 1-q^{2d_{32}}} }
 \cr
0 &0 &0 &0 &\displaystyle{ {q-q^{2d_{23}-1}\over 1-q^{2d_{23}}} }
&\displaystyle{ {q-q^{-1}\over 1-q^{2d_{32}}} } 
 \cr
}
$$
and
$$
\phi_S(T_{i+1})=\pmatrix{
\displaystyle{ {q-q^{-1}\over 1-q^{2d_{23}}} } 
&0 &\displaystyle{ {q-q^{2d_{32}-1}\over 1-q^{2d_{32}}} }
&0 &0 &0 \cr
0 &\displaystyle{ {q-q^{-1}\over 1-q^{2d_{13}}} } 
&0 &0 &\displaystyle{ {q-q^{2d_{31}-1}\over 1-q^{2d_{31}}} }
&0\cr
\displaystyle{ {q-q^{2d_{23}-1}\over 1-q^{2d_{23}}} }
&0 &\displaystyle{ {q-q^{-1}\over 1-q^{2d_{32}}} } 
&0 &0 &0 \cr
0 &0 &0 &\displaystyle{ {q-q^{-1}\over 1-q^{2d_{12}}} } 
&0 &\displaystyle{ {q-q^{2d_{21}-1}\over 1-q^{2d_{21}}} }
 \cr
0 &\displaystyle{ {q-q^{2d_{13}-1}\over 1-q^{2d_{13}}} }
&0 &0 &\displaystyle{ {q-q^{-1}\over 1-q^{2d_{31}}} } 
&0 \cr
0 &0 &0 &\displaystyle{ {q-q^{2d_{12}-1}\over 1-q^{2d_{12}}} }
&0 &\displaystyle{ {q-q^{-1}\over 1-q^{2d_{21}}} } 
\cr
}
$$
where $d_{k\ell}=c(L(i+k-1))-c(L(i+l-1))$.
In each case we compute directly the products
$\phi_S(T_i)\phi_S(T_{i+1})\phi_S(T_i)$
and 
$\phi_S(T_{i+1})\phi_S(T_i)\phi_S(T_{i+1})$
and verify that they are equal.   (This proof of the braid relation is, in all
essential aspects, the same as that used by Hoefsmit [Ho], Wenzl [Wz] and Ariki
and Koike [AK].  For a more elegant but less straightforward proof of this
relation see the proof of Theorem 3.1 in [Ra3].)
\endpf

\noindent
{\smallcaps step 2.}  The module $\tilde H^{(c,\lambda/\mu)}$
is irreducible.
\pf
Let $L$ be a standard tableaux of shape $\lambda/\mu$
and define 
$$\pi_L = \prod_{i=1}^n \prod_{P\ne L} 
{x_i-q^{2c(P(i))}\over q^{2c(L(i))}-q^{2c(P(i))}},$$
where the second product is over all standard tableaux $P$ of shape
$\lambda/\mu$ which are not equal to $L$.
Then $\pi_L$ is an element of $\tilde H_n$ such that
$$\pi_L v_Q = \delta_{LQ}v_L,$$
for all standard tableaux $Q$ of shape $\lambda/\mu$.
This follows from the formula for the 
action of $x_i$ on $\tilde H^{(c,\lambda/\mu)}$
and the fact that the sequence
$(q^{2c(L(1))}, \ldots, q^{2c(L(n))})$ completely determines the standard
tableau $L$ (Lemma 2.2). 

Let $N$ be a nonzero submodule of $\tilde H^{(c,\lambda/\mu)}$ and
let $v=\sum_{Q} a_Qv_Q$ be a nonzero element of $N$. 
Let $L$ be a standard tableau such that the coefficient
$a_L$ is nonzero. Then $\pi_Lv = a_Lv_L$ and so $v_L\in N$.

By Proposition 2.1 we may identify the set $\cF^{\lambda/\mu}$
with an interval in $S_n$ (under Bruhat order).  Under this identification
the minimal element is the column reading tableau $C$
and there is a chain 
$C<s_{i_1}C<\cdots<s_{i_p}\cdots s_{i_1}C=L$ such that all elements
of the chain are standard tableaux of shape $\lambda/\mu$.  Then,
by the definition of the $\tau_i$-operators,
$$\tau_{i_1}\cdots \tau_{i_p}v_L = \kappa v_C,$$
for some constant $\kappa\in \CC^*$. It follows that
$v_C\in N$.

Let $Q$ be an arbitrary standard tableau of shape $\lambda/\mu$.
Again, there is a chain 
$C<s_{j_1}C<\cdots<s_{j_p}\cdots s_{j_1}C=Q$ of standard tableaux
in $\cF^{\lambda/\mu}$ and we have
$$\tau_{j_p}\cdots \tau_{j_1}v_C = \kappa'v_Q,$$
for some $\kappa'\in \CC^*$.  Thus $v_Q\in N$.  It follows
that $N=\tilde H^{(c,\lambda/\mu)}$.    Thus 
$\tilde H^{(c,\lambda/\mu)}$ is irreducible.
\endpf

\noindent
{\smallcaps step 3.}  The modules $\tilde H^{(c,\lambda/\mu)}$ are
nonisomorphic.
\pf  Each of the modules $\tilde H^{(c,\lambda/\mu)}$ has a unique basis
(up to multiplication of each basis vector by a constant) of simultaneous
eigenvectors for the $x_i$.  Each basis vector is determined by
its weight, the sequence of eigenvalues $(t_1,\ldots, t_n)$ given by
$$x_i v_t = t_i v_t,
\qquad\hbox{for $1\le i\le n$.}$$
By the definition of the action of the $x_i$,
a weight of $\tilde H^{(c,\lambda/\mu)}$ is equal to
$(q^{2c(L(1))},\ldots,q^{2c(L(n))})$ 
for some standard tableau $L$.  By Lemma 2.2, both the
standard tableau $L$ and the placed skew shape $(c,\lambda/\mu)$ 
are determined uniquely by this weight.
Thus no two of the modules $\tilde H^{(c,\lambda/\mu)}$ can be isomorphic.
\endpf

\medskip\noindent
{\smallcaps step 4.}  If $t=(t_1,\ldots, t_n)$ is the weight of a calibrated 
$\tilde H_n$-module $M$ then  $t=(q^{2c(L(1))},\ldots,q^{2c(L(n))})$ for some
standard tableau $L$ of placed skew shape.  

\pf
Let $m_t$ be a weight vector in $M$ of weight $t=(t_1,\ldots, t_n)$, i.e.
$$x_im_t = t_im_t, \qquad\hbox{for all $1\le i\le n$.}$$ 
We want $L$ such that $(q^{2c(L(1))},\ldots,q^{2c(L(n))})=(t_1,\ldots, t_n)$.
We shall show that if
$t_i=t_j$ for $i<j$ then there exist $k$ and
$\ell$ such that $i<k<\ell<j$, $t_k = q^{\pm2}t_i$ and $t_\ell=q^{\mp2}t_i$.
This will show that if there are two adjacent boxes of $L$ in the
same diagonal then these boxes must be contained in a complete $2\times 2$ block,
i.e. if there is a confuguration in $L$ of the form
$$\beginpicture
\setcoordinatesystem units <0.5cm,0.5cm>         
\setplotarea x from 0 to 2, y from 0 to 2    
\linethickness=0.5pt                          
\put{$i$} at .5 1.5
\put{$j$} at 1.5 .5
\putrule from 0 1 to 2 1          %
\putrule from 0 2 to 1 2          
\putrule from 1 0 to 2 0          
\putrule from 0 1 to 0 2        %
\putrule from 1 0 to 1 2        %
\putrule from 2 0 to 2 1        
\endpicture  
\qquad\hbox{then $L$ must contain}\qquad
\beginpicture
\setcoordinatesystem units <0.5cm,0.5cm>         
\setplotarea x from 0 to 2, y from 0 to 2    
\linethickness=0.5pt                          
\put{$i$} at .5 1.5
\put{$j$} at 1.5 .5
\put{$k$} at 1.5 1.5
\put{$\ell$} at .5 .5
\putrule from 0 1 to 2 1          %
\putrule from 0 2 to 2 2          
\putrule from 0 0 to 2 0          
\putrule from 0 0 to 0 2        %
\putrule from 1 0 to 1 2        %
\putrule from 2 0 to 2 2        
\endpicture  
\qquad\hbox{or}\qquad
\beginpicture
\setcoordinatesystem units <0.5cm,0.5cm>         
\setplotarea x from 0 to 3, y from 0 to 2    
\linethickness=0.5pt                          
\put{$i$} at .5 1.5
\put{$j$} at 1.5 .5
\put{$\ell$} at 1.5 1.5
\put{$k$} at .5 .5
\putrule from 0 1 to 2 1          %
\putrule from 0 2 to 2 2          
\putrule from 0 0 to 2 0          
\putrule from 0 0 to 0 2        %
\putrule from 1 0 to 1 2        %
\putrule from 2 0 to 2 2        
\endpicture  
.
$$
This is sufficient to guarantee that $L$ is of skew shape.

Let $j>i$ be such that $t_j=t_i$ and $j-i$ is minimal.  The argument is by
induction on the value of $j-i$.

\smallskip\noindent
{\it Case 1: $j-i=1$.}  Then $m_t$ and
$T_im_t$ are linearly independent. If they were not we would have
$T_im_t = am_t$ which would give
$$t_i am_t = x_iT_im_t
=(T_ix_{i+1} - (q-q^{-1})x_{i+1}) m_t
=(at_{i+1} - (q-q^{-1})t_{i+1})m_t = (a-(q-q^{-1}))t_im_t.$$
Since $q-q^{-1}\ne 0$, this equation implies that $t_i = 0$ which is a
contradiction.   Now the relations (1.1d) and (1.4) show that
$$\eqalign{
x_iT_im_t &=t_i(T_im_t - (q-q^{-1})m_t), \cr
x_{i+1}T_i m_t &= t_{i+1}( T_im_t+(q-q^{-1}) m_t), \cr
x_jT_im_t &= t_j T_im_t, \qquad\hbox{for all $j\ne i,i+1$.} \cr
}$$
It follows that $T_im_t$ is an element of $M_t^{\rm gen}$ but not an element of
$M_t$.  This is a contradiction to the fact that $M$
is calibrated.  So this case is not possible, i.e. $t_{i+1}$ {\it cannot} equal
$t_i$.

\smallskip\noindent
{\it Case 2: $j-i=2$.}
Since $t_i\ne t_{i+1}$ and $m_t$ is a weight vector,
the vector
$$m_{s_it} = \left(T_i - {(q-q^{-1})t_{i+1}\over t_{i+1}-t_i}\right)m_t$$
is a weight vector of weight $t'=s_it$ (see (3.7a)). Then
$t'_i=t'_{i+1}$ and so, by Case 1, $m_{s_it}=0$.  
This implies that 
$$T_im_t = {(q-q^{-1})t_{i+1}\over t_{i+1}-t_i}m_t.$$
By equation (1.3), all eigenvalues of $T_i$ are either $q$ or $-q^{-1}$.
Thus $T_im_t = \pm q^{\pm 1}m_t$ and so $t_i=q^{\pm2}t_{i+1}$.
A similar argument shows that
$$m_{s_{i+1}t}=
\left(T_{i+1} - {(q-q^{-1})t_{i+2}\over t_{i+2}-t_{i+1}}\right)m_t$$ 
must be $0$ and thus that   
$$
T_{i+1}m_t = {(q-q^{-1})t_{i+2}\over t_{i+2}-t_{i+1}}m_t 
= {(q-q^{-1}) t_i\over t_i-t_{i+1}} m_t 
= \mp q^{\mp1}m_t.$$
From $T_im_t=\pm q^{\pm1}m_t$ and $T_{i+1}m_t=\mp q^{\mp1}m_t$ we get
$$\pm q^{\pm 1}m_t = T_iT_{i+1}T_im_t = T_{i+1}T_iT_{i+1}m_t 
= \mp q^{\mp1}m_t.$$
This is impossible since $q$ is not a root of unity.  So this case is not
possible, i.e. $t_{i+2}$ {\it cannot} equal $t_i$.

\smallskip\noindent
{\it Induction step.}   Assume that $i$ and $j$ are such that $t_i=t_j$
and the value $j-i$ is minimal such that this is true.

If $t_{j-1}\ne q^{\pm2}t_j$ then the vector
$$m_{s_j t} = \left(T_j - {(q-q^{-1})t_j\over t_{j-1}-t_j}\right) m_t$$
is a weight vector of weight $t'=s_jt$ and by (3.7b) this vector is
nonzero.  Since $t'_i=t_i=t_j=t'_{j-1}$ we can apply
the induction hypothesis to conclude that there are $k$ and $\ell$
with $i<k<\ell<j-1$ such that $t'_k=q^{\pm2}t'_i$ and $t'_\ell =
q^{\mp2}t'_i$.  This implies that $t_k=q^{\pm2}t_i$ and
$t_\ell=q^{\mp2}t_i$.

Similarly, if $t_i \ne q^{\pm2}t_{i+1}$
then the vector
$$m_{s_i t} = \left(T_i - {(q-q^{-1})t_{i+1}\over t_{i+1}-t_i}\right) m_t$$
is a weight vector of weight $t'=s_it$ and by (3.7b) this vector is
nonzero.  Since $t'_{i+1}=t_i=t_j=t'_j$ we can apply
the induction hypothesis to conclude that there are $k$ and $\ell$
with $i+1<k<\ell<j$ such that $t'_k=q^{\pm2}t'_j$ and $t'_\ell =
q^{\mp2}t'_j$.  This implies that $t_k=q^{\pm2}t_i$ and
$t_\ell=q^{\mp2}t_i$.

If we are not in either of the previous cases then
$t_{i+1}=q^2t_i$ or $t_{i+1}=q^{-2}t_i$ and $t_{j-1}=q^2t_j$ or
$t_{j-1}=q^{-2}t_j$.   We cannot have $t_{i+1}=t_{j-1}$ since the $i$ and $j$
are such that $j-i$ is minimal such that $t_i=t_j$.
Thus $q^{\pm2}t_{i+1}=q^{\mp 2}t_{j-1}=t_i.$
\endpf

\medskip\noindent
{\smallcaps step 5.}   Suppose that $M$ is an irreducible calibrated $\tilde
H_n$-module and that $m_t$ is a weight vector in $M$ with weight $t=(t_1,\ldots,
t_n)$ such that $t_i=q^{\pm 2}t_{i+1}$.  Then $\tau_i m_t = 0$.

\pf  Assume that $m_{s_it}=\tau_im_t\ne 0$.  Then, by the second identity in
Proposition 3.2 (b),
$\tau_i m_{s_it} = 0$.  Since $M$ is
irreducible there must be some sequence of $\tau$-operators such that
$$\tau_{i_1}\ldots \tau_{i_p} m_{s_it} = \kappa m_t,$$
with $\kappa\in \CC^*$.   Assume that $\tau_{i_1}\cdots \tau_{i_p}$ is
a minimal length sequence such that this is true.  We have $s_{i_1}\cdots
s_{i_p}s_i t=t$.   

Assume that $s_{i_1}\cdots s_{i_p}s_i\ne 1$.  Then there must be
$1\le i<j\le n$ such that $t_i=t_j$ (because some nontrivial permutation
of the $t_i$ fixes $t$).  Since $s_{i_1}\cdots s_{i_p}s_i$ is of minimal
length such that it fixes $t$ it must be a transposition $(i,j)$
for some $i<j$ such that $t_i=t_j$.  Furthermore there does not exist
$i<k<j$ such that $t_i=t_k$.  The element $s_{i_1}\cdots s_{i_p}s_i$
switches the $t_i$ and the $t_j$ in $t$.  In the process of doing this switch
by a sequence of simple transpositions there must be some point where $t_i$ and
$t_j$ are adjacent and thus there must be some $\ell$ such that
$s_{i_\ell}(s_{i_{\ell+1}}\cdots s_{i_p}s_it)=s_{i_{\ell+1}}\cdots s_{i_p}s_it$.
Then 
$$m_{t'}=\tau_{i_{\ell+1}}\cdots \tau_{i_\ell}\tau_im_t$$
is a nonzero weight vector in $M$ of weight $t'=s_{i_{\ell+1}}\cdots
s_{i_p}s_it$. Since $s_{i_\ell}t'=t'$ it follows that
$t_{i_\ell}=t_{i_\ell+1}$.   Since $M$ is calibrated this is a contradiction to
(Case 1 of) {\smallcaps step 4}.  

So $s_{i_1}\cdots s_{i_p}s_i=1$.  Let $k$ be minimal such that
$s_{i_1}\cdots s_{i_k}$ is not reduced.  Assume $p>1$. Then we can use the braid
relations (the third and fourth lines of Proposition 3.2 (b)) to write 
$$\kappa m_t = \tau_{j_1}\cdots
\tau_{j_{k-2}}\tau_{i_k}\tau_{i_k}\tau_{i_{k+1}}\cdots\tau_{i_p}\tau_im_t,$$
for some $j_1,\ldots, j_{k-2}$.
Then, by the second line of Proposition 3.2 (b), 
$$\kappa'm_t = \tau_{j_1}\cdots
\tau_{j_{k-2}}\tau_{i_{k+1}}\cdots\tau_{i_p}\tau_im_t,$$
for some $\kappa'\in \CC^*$.  This is a contradiction to the minimality
of the length of the sequence $\tau_{i_1}\cdots \tau_{i_p}$.

So $p=1$, $i_p=i$ and $\tau_i\tau_i m_t = \kappa m_t$.  This is a contradiction
since the second identity in Proposition 3.2 (b) and the assumption that
$t_i=q^{\pm 2}t_{i+1}$ imply that
$\tau_i\tau_i=0$.  So $\tau_im_t=0$.
\endpf

\medskip\noindent
{\smallcaps step 6.}  An irreducible calibrated $\tilde H_n$-module $M$ is
isomorphic to $\tilde H^{(c,\lambda/\mu)}$ for some placed skew shape
$(c,\lambda/\mu)$.

\pf
Let $m_t$ be a nonzero weight vector in $M$.  Since $M$ is calibrated
{\smallcaps step 4} implies that there is a placed skew shape $(c,\lambda/\mu)$
and a standard tableau $L$ of shape $\lambda/\mu$ such that
$t=(q^{2c(L(1))},\ldots,q^{2c(L(n))})$.  Let us write $m_L$ in place of $m_t$.
Let $C$ be
the column reading tableau of shape $\lambda/\mu$.  It follows from Proposition
2.1 that there is a chain
$C, s_{i_1}C,\ldots, s_{i_p}\cdots s_{i_1}C=L$ of standard tableaux of shape
$\lambda/\mu$.  By (3.7b), all of the $\tau_{i_j}$
in this sequence are bijections and so 
$$m_C = \tau_{i_1}\cdots \tau_{i_p} m_L$$ 
is a nonzero weight vector in $M$.  Similarly, if $Q$ is any other standard
tableau of shape $\lambda/\mu$ then there is a chain
$C, s_{j_1}C,\ldots, s_{j_p}\cdots s_{j_1}C=Q$ and so
$$m_Q = \tau_{j_p}\cdots \tau_{j_1} m_C$$ 
is a nonzero weight vector in $M$.
Finally, by {\smallcaps step 5},  $\tau_i m_Q= 0$ if $s_iQ$ is not standard 
(since $q^{2c(Q(i))}=q^{\pm 2}q^{2c(Q(i+1))}$) and so the span of the
vectors
$\{ m_Q\ |\
\hbox{$Q$ a standard tableau of shape $\lambda/\mu$}\}$ is a submodule of $M$. 
Since $M$ is irreducible this must be all of $M$ and the map
$$\matrix{
M &\longrightarrow &\tilde H^{(c,\lambda/\mu)} \cr
\tau_v m_{C} &\longmapsto &\tau_v v_{C} \cr
}$$
is an isomorphism of $\tilde H_n$-modules.
\endpf

\noindent
This completes the proof of Theorem 4.1.

\bigskip

\section 5. ``Garnir relations'' and an analogue of Young's natural basis

\bigskip
Each of the modules $\tilde H^{(c,\lambda/\mu)}$ constructed in Theorem
4.1 has two natural bases:
\item{(a)}  The ``seminormal basis'' $\{v_L\ |\ L\in\cF^{\lambda/\mu}\}$
which, up to multiplication of each basis vector by a constant, is given by
$$\tilde v_L = \tau_w v_C,\formula$$
where $C$ is the column reading tableau of shape $\lambda/\mu$, $w$ is the
permutation such that $L=wC$ and $\tau_w$ is as in (3.3).
\item{(b)}  The ``natural basis'' $\{ n_L \ |\ L\in \cF^{\lambda/\mu} \}$
given by 
$$n_L = T_w v_C, \formula$$
where $C$ is the column reading tableau of shape $\lambda/\mu$, $w$ is the
permutation such that $L=wC$ and $T_w$ is as defined in (1.6).

\prop  Let $(c,\lambda/\mu)$ be a placed skew shape with $n$ boxes and let
$\tilde H^{(c,\lambda/\mu)}$ be the $\tilde H_n$-module defined in Theorem 4.1.
Let $C$ be the  column reading tableau of shape $\lambda/\mu$ and let $wC$
denote the tableau $C$ with the entries permuted according to the permutation
$w$. For each standard tableau $L$ let $n_L$ be defined by formula (5.2).
Then $\{ n_L\ |\  L\in \cF^{\lambda/\mu}\}$ is a basis of 
$\tilde H^{(c,\lambda/\mu)}$.
\pf
If $L$ is a standard tableau of shape $\lambda/\mu$ let $\tilde v_L$ be as given
in (5.1).  It follows from the formulas defining the module
$\tilde H^{(c,\lambda/\mu)}$ that the basis $\{\tilde v_L\ |\ L\in
\cF^{\lambda/\mu}\}$ is simply a renormalized version of the basis 
$\{v_L\ |\ L\in \cF^{\lambda/\mu}\}$, i.e there are constants
$\kappa_L\in \CC^*$ such that $\tilde v_L = \kappa_Lv_L$.

Let $s_{i_1}\cdots s_{i_p}=w$ be a reduced
word for $w$.  Then, with notations as in Theorem 4.1,
$$
\tilde v_L = \tau_{i_1}\ldots \tau_{i_p} v_C  
= (T_{i_1} - (T_{i_1})_{L_2L_2})(T_{i_2} - (T_{i_1})_{L_3L_3})
\cdots
(T_{i_p}-(T_{i_p})_{L_pL_p})n_C,
$$
where $L_j = s_{i_{j+1}}\cdots s_{i_p}C$.  Expanding this expression yields
$$\tilde v_L=(T_w+\sum_{u<w} b_u T_u)n_C
=n_L + \sum_{u<w} b_u n_{uC},$$
for some constants $b_u\in \CC$.  The second equality is a consequence of the
fact that, by Proposition 2.1, the tableaux
$uC$, $u<w$, are always standard.  Since $\{\tilde v_L\ |\ L\in
\cF^{\lambda/\mu}\}$ is a basis it follows from the triangular relation above
that $\{n_L\ |\ L\in \cF^{\lambda/\mu}\}$ is also a basis of $\tilde
H^{(c,\lambda/\mu)}$.
\endpf

The construction of $\tilde H^{(c,\lambda/\mu)}$ in Theorem 4.1 makes the
notational assumption that $v_L=0$, whenever $L$ is not a standard tableau.
Formula (5.1) can be used as a definition of $\tilde v_L$ even in the case
when $L$ is not standard and, from the definition of the action in Theorem 4.1,
$$\tilde v_L = 0,\qquad\hbox{if $L$ is not standard.}\formula$$
Theorem 5.5 proves that these relations, when expanded in terms of the
basis $\{n_L \ |\ L\in \cF^{\lambda/\mu}\}$ are exactly the classical Garnir
relations!! 

Let $\lambda/\mu$ be a skew shape.  A pair of adjacent boxes in the same row
of $\lambda/\mu$ determines a {\it snake} in $\lambda/\mu$ consisting
of the boxes in the pair, all the boxes above the righthand box of this pair, and
all the boxes below the lefthand box of the pair.  See the picture in 
Theorem 5.5 (b).

\thm  (``Garnir relations'') \ \  Let $(c,\lambda/\mu)$ be a placed skew shape
and let $\tilde H^{(c,\lambda/\mu)}$ be the $\tilde H_n$-module defined in
Theorem 4.1. Let $\{n_L\ |\ L\in \cF^{\lambda/\mu}\}$ be the basis of $\tilde
H^{(c,\lambda/\mu)}$ defined by Proposition 5.3 and let $C$ be the column
reading tableau of shape $\lambda/\mu$.  
\item{(a)}  If $i$ and $i+1$ are entries in the same column of $C$
then
$$T_i\; n_C = -q^{-1} n_C.$$
\item{(b)}  Fix a snake in $\lambda/\mu$.  Let $P$ be the standard tableau which
has all entries the same as $C$ except that the entries in the snake
are entered in row reading order instead of in column reading order.
$$
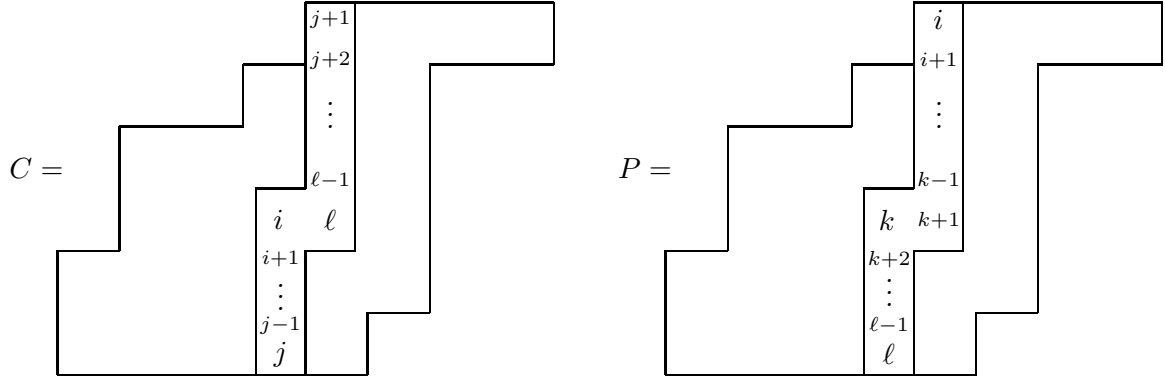
\beginpicture
\setcoordinatesystem units <0.5cm,0.5cm>         
\setplotarea x from 0 to 14, y from 0 to 9    
\put{$C=$} at 1 5
\put{$j$}  at 6.9 0.5
\put{$\scriptstyle{j-1}$} at 6.9 1.2
\put{$\vdots$} at 6.9 2.1
\put{$\scriptstyle{i+1}$} at 6.9 2.8
\put{$i$}  at 6.85 3.75
\put{$\ell$}  at 8.1 3.75
\put{$\scriptstyle{\ell-1}$}  at 8.1 4.7
\put{$\vdots$} at 8.1 6.5
\put{$\scriptstyle{j+2}$}  at 8.1 7.6
\put{$\scriptstyle{j+1}$}  at 8.1 8.6
\putrule from 7.5 9 to 13.5 9               %
\putrule from 6 7.5 to 7.5 7.5              %
\putrule from 10.5 7.5 to 13.5 7.5          
\putrule from 3 6 to 6 6                    %
\putrule from 1.5 3 to 3 3                  %
\putrule from 9 1.5 to 10.5 1.5             %
\putrule from 1.5 0 to 9 0                  %
\putrule from 1.5 0 to 1.5 3            %
\putrule from 3 3 to 3 6                
\putrule from 6 6 to 6 7.5              %
\putrule from 7.5 7.5 to 7.5 9          %
\putrule from 9 0  to  9 1.5            %
\putrule from 10.5 1.5 to 10.5 7.5      %
\putrule from 13.5 7.5 to 13.5 9        %
\putrule from 6.3 0 to 6.3 4.5  
\putrule from 7.5 0 to 7.5 3
\putrule from 6.3 4.5 to 7.5 4.5
\putrule from 7.5 4.5 to 7.5 7.5
\putrule from 7.5 3 to 8.7 3
\putrule from 8.7 3 to 8.7 9 
\endpicture
\quad
\beginpicture
\setcoordinatesystem units <0.5cm,0.5cm>         
\setplotarea x from 0 to 14, y from 0 to 9    
\put{$P=$} at 1 5
\put{$\ell$}  at 6.9 0.5
\put{$\scriptstyle{\ell-1}$} at 6.9 1.2
\put{$\vdots$} at 6.9 2.2
\put{$\scriptstyle{k+2}$} at 6.9 2.8
\put{$k$}  at 6.85 3.75
\put{$\scriptstyle{k+1}$}  at 8.1 3.75
\put{$\scriptstyle{k-1}$}  at 8.1 4.7
\put{$\vdots$} at 8.1 6.5
\put{$\scriptstyle{i+1}$}  at 8.1 7.6
\put{$i$}  at 8.1 8.6
\putrule from 7.5 9 to 13.5 9               %
\putrule from 6 7.5 to 7.5 7.5              %
\putrule from 10.5 7.5 to 13.5 7.5          
\putrule from 3 6 to 6 6                    %
\putrule from 1.5 3 to 3 3                  %
\putrule from 9 1.5 to 10.5 1.5             %
\putrule from 1.5 0 to 9 0                  %
\putrule from 1.5 0 to 1.5 3            %
\putrule from 3 3 to 3 6                
\putrule from 6 6 to 6 7.5              %
\putrule from 7.5 7.5 to 7.5 9          %
\putrule from 9 0  to  9 1.5            %
\putrule from 10.5 1.5 to 10.5 7.5      %
\putrule from 13.5 7.5 to 13.5 9        %
\putrule from 6.3 0 to 6.3 4.5  
\putrule from 7.5 0 to 7.5 3
\putrule from 6.3 4.5 to 7.5 4.5
\putrule from 7.5 4.5 to 7.5 7.5
\putrule from 7.5 3 to 8.7 3
\putrule from 8.7 3 to 8.7 9 
\endpicture
$$
Let $S_A$, $S_B$ and $S_{A\cup B}$ be the subgroups of $S_n$
consisting of the
permutations of
$$A=\{i,i+1,\ldots,j\}, \qquad
B=\{j+1,\ldots, \ell-1,\ell\}
$$
and $A\cup B$, respectively, and let $S_{A\cup B}/(S_A\times S_B)$ be the set of
minimal length coset representatives of cosets of $S_A\times S_B$ in
$S_{A\cup B}$.  The elements of $S_{A\cup B}/(S_A\times S_B)$ are sometimes
called the ``shuffles'' of $A$ and $B$.
Then
$$\eqalign{
0&=
\left(\sum_{u\in S_{A\cup B}/(S_A\times S_B)}
(-q)^{\ell(w_{A\cup B})-\ell(w_Aw_B)-\ell(u)} T_u\right)
n_C \cr
&=T_k n_P + \sum_{u\le P}
(-q)^{\ell(w_{A\cup B})-\ell(w_Aw_B)-\ell(u)} n_{uC},
\cr}
$$
where $\ell(w_{A\cup B}) = {\ell-i+1\choose 2}$, $\ell(w_Aw_B)={j-i+1\choose
2}{\ell-j\choose 2}$ and the last sum is over all standard tableaux $uC$ which
are obtained from $C$ by permuting entries which are in the snake.
\pf
Part (a) follows immediately from the definition of $\tilde H^{(c,\lambda/\mu)}$
in Theorem 4.1.
\smallskip\noindent
(b)  
The subgroups $S_A$, $S_B$ and $S_{A\cup B}$ have longest elements
$$
w_A = \pmatrix{
i &i+1 &\cdots &j-1 &j \cr
j &j-1 &\cdots &i+1 &i \cr}, \qquad\qquad
w_B = \pmatrix{
j+1 &j+2 &\cdots &\ell-1 &\ell \cr
\ell &\ell-1 &\cdots &j+2 &j+1 \cr},
$$
$$
w_{A\cup B} = \pmatrix{
i &i+1 &\cdots &\ell-1 &\ell \cr
\ell &\ell-1 &\cdots &i+1 &i \cr},
$$
with lengths $\ell(w_A) = {j-i+1\choose 2}$, $\ell(w_B)={\ell-j\choose
2}$ and $\ell(w_{A\cup B})={\ell-i+1\choose 2}$, respectively, and
$w_Aw_B$ is the longest element of $S_A\times S_B\subseteq S_{A\cup B}$.
Let $t=(t_1,\ldots, t_n)=(a, aq^2, aq^4,\ldots, aq^{2(n-1)})$
where $a=q^{2(c(C(j))-(i-1))}\in \CC^*$. The positions of $C(i), \ldots C(\ell)$
in $\lambda/\mu$ are such that
$$\eqalign{
(q^{2c(C(i))}, &\ldots, q^{2c(C(\ell))}) \cr
&=w_Aw_B(q^{2c(C(j))},q^{2c(C(j-1))},\ldots, q^{2c(C(i))},
q^{2c(C(\ell))}, q^{2c(C(\ell-1))},\ldots, q^{2c(C(j+1))}) \cr
&=w_Aw_B(q^{2c(C(j))},q^{2c(C(j))}q^2,\ldots, q^{2c(C(j))}q^{2(\ell-i+1)})
\cr
&=w_Aw_B(t_i,\ldots,t_\ell). \cr}$$
By Proposition 3.6 and the fact that $T_s n_C = -q^{-1}n_C$ for all $i\le s\le
\ell$, $s\ne j$,
$$
\tau_{w_Aw_B}\big\vert_{t} n_C
= \left(
\sum_{w\in S_A\times S_B} (-q)^{\ell(w_Aw_B)-\ell(w)} T_w \right)
n_C
= [j-i+1]![\ell-j]! n_C,
$$
where $[r]! =[r][r-1]\cdots [2][1]$ and $[p] = (q^p-q^{-p})/(q-q^{-1})$.

Let $\pi$ be the permutation such that $P=\pi C$ and let $\tilde v_P = \tau_\pi
v_C$.  Then
$s_i\pi w_Aw_B = w_{A\cup B}$ and, by (3.5),
$$\eqalign{
\tau_i \tilde v_P  &= \tau_i\tau_\pi v_C
= \tau_i\tau_\pi\big\vert_{w_Aw_B t} n_C \cr
&= {1\over [j-i+1]![\ell-j]!} \tau_i\tau_\pi\big\vert_{w_Aw_B t}
\tau_{w_Aw_B}\big\vert_{t} n_C \cr
&= {1\over [j-i+1]![\ell-j]!} \tau_i\tau_\pi
\tau_{w_Aw_B}\big\vert_{t} n_C \cr
&= {1\over [j-i+1]![\ell-j]!}
\tau_{w_{A\cup B}}\big\vert_{t} n_C \cr
&= {1\over [j-i+1]![\ell-j]!}
\left(
\sum_{w\in S_{A\cup B}} (-q)^{\ell(w_{A\cup B})-\ell(w)}
T_w \right) n_C.
\cr}$$
Each element $w\in S_{A\cup B}$ has a unique expression
$w=uv$ such that $v\in S_A\times S_B$ and $\ell(w)=
\ell(u)+\ell(v)$.  The left factor $u$ is the minimal length
representative of the coset $w(S_A\times S_B)$ in $S_{A\cup B}$.
Then
$$\eqalign{
\sum_{w\in S_{A\cup B}} (-q&)^{\ell(w_{A\cup B})-\ell(w)} T_w \cr
&=
\left(\sum_{u\in S_{A\cup B}/(S_A\times S_B)}
(-q)^{\ell(w_{A\cup B})-\ell(w_Aw_B)-\ell(u)} T_u\right)
\left(\sum_{v\in S_A\times S_B}
(-q)^{\ell(w_Aw_B)-\ell(v)} T_v
\right). \cr}
$$
It follows from the formulas for the action on $\tilde H^{(c,\lambda/\mu)}$
that the element $\tilde v_P = \tau_\pi v_C$ is a nonzero multiple of the basis
element $v_P$.  Furthermore, by (5.4), $\tau_k\tilde v_P=0$.  So
$$\eqalign{
0&=\tau_i \tilde v_P  \cr
&= {1\over [j-i+1]![\ell-j]!}
\left(\sum_{u\in S_{A\cup B}/(S_A\times S_B)}
(-q)^{\ell(w_{A\cup B})-\ell(w_Aw_B)-\ell(u)} T_u\right) \cr
&\qquad\qquad\qquad\qquad\qquad\qquad\qquad\qquad\qquad \times
\left(\sum_{v\in S_A\times S_B} (-q)^{\ell(w_Aw_B)-\ell(v)} T_v
\right)n_C \cr
&=
\left(\sum_{u\in S_{A\cup B}/(S_A\times S_B)}
(-q)^{\ell(w_{A\cup B})-\ell(w_Aw_B)-\ell(u)} T_u\right)
n_C \cr
}$$
For each $u\in S_{A\cup B}/(S_A\times S_B)$ except $u=w_{A\cup B}w_Aw_B$, the
tableau $uC$ is standard.  In fact these are exactly the standard tableaux $Q$
which are obtained by permuting entries of $C$ which are in the snake.
The tableau $w_{A\cup B}w_Aw_BC = s_kP$.  Note that $\ell(w_{A\cup
B}w_Aw_B)=\ell(w_{A\cup B})-\ell(w_Aw_B)$.  Thus we have
$$0=T_k n_P + \sum_{u\le P}
(-q)^{\ell(w_{A\cup B})-\ell(w_Aw_B)-\ell(u)} n_{uC},
$$
where the sum over all standard tableaux $uC$ which are obtained from $C$ by
permuting entries which are in the snake.
\endpf

\prop  Let $(c,\lambda/\mu)$ be a placed skew shape and let
$\{n_L \ |\ \hbox{$L$ is a standard tableau of}$ $\hbox{ shape $\lambda/\mu$}\}$
be the basis of the $\tilde H_n$-module $\tilde H^{(c,\lambda/\mu)}$ which is
defined by Proposition 5.3.  
\item{(a)}  If
$w\in S_n$ and
$L$ is a standard tableau of shape $\lambda/\mu$ then
$$T_w n_L = \sum_{Q} b_Q n_Q,\qquad
\hbox{with coefficients $b_Q\in \ZZ[q,q^{-1}]$.}$$
\item{(b)}  Assume that the content function $c$ takes values in $\ZZ$.  If
$1\le i\le n$ and $L$ is a standard tableau of shape $\lambda/\mu$ then
$$x_i n_L = \sum_{Q} b'_Q n_Q,\qquad
\hbox{with coefficients $b'_Q\in \ZZ[q,q^{-1}]$.}$$

\pf
(a)  It is sufficient to show that for all $1\le d\le n$ and all standard
tableaux $L$ of shape $\lambda/\mu$ we have 
$$T_d n_L =\sum_Q b_Qn_Q,$$
with coefficients $b_Q\in \ZZ[q,q^{-1}]$.  For notational convenience
let us identify each standard tableau $L$ of shape $\lambda/\mu$
with the permutation $w\in S_n$ such that $wC=L$, where $C$ is the column reading
tableau of shape $\lambda/\mu$.  

From the definitions
$$\matrix{
T_d n_L = n_{s_dL}, \hfill
\qquad &\hbox{if $\ell(s_dL)=\ell(L)+1$ and $s_dL$ is standard,}\hfill \cr
T_dn_L = (q-q^{-1})n_L + n_{s_dL}, \hfill
\qquad &\hbox{if $\ell(s_dL)=\ell(L)-1$.}\hfill \cr
}$$
Let $L$ be a standard tableau such that $\ell(s_dL) = \ell(L) + 1$ and 
$s_dL$ is not standard.  The pair of boxes where the nonstandardness
in $s_dL$ occurs in either a row or a column:
$$
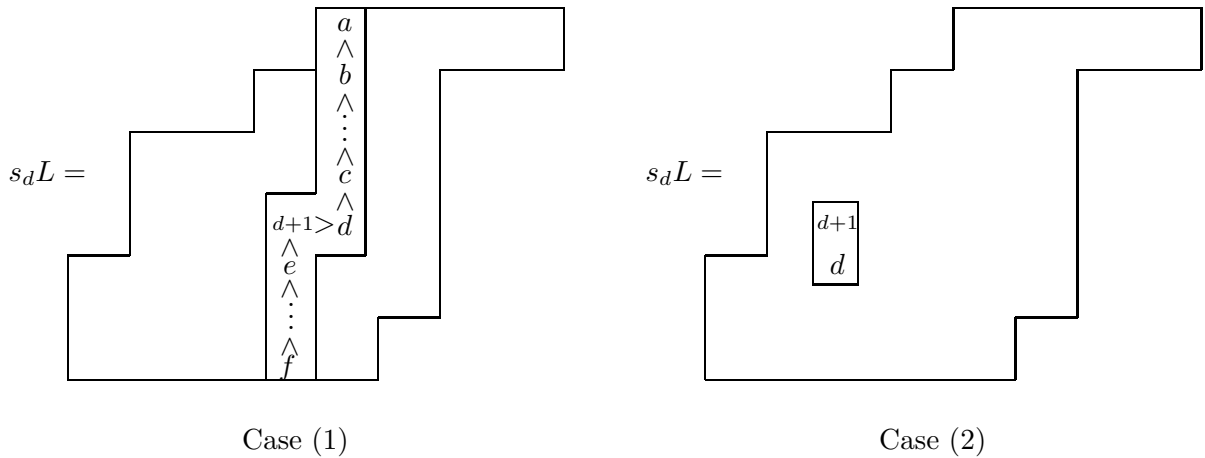
\matrix{
\beginpicture
\setcoordinatesystem units <0.5cm,0.5cm>         
\setplotarea x from 0 to 14, y from 0 to 9    
\put{$s_dL=$} at 1 5
\put{$f$}  at 6.8 0.3
\put{$\wedge$} at 6.87 .8
\put{$\vdots$} at 6.9 1.75
\put{$\wedge$} at 6.87 2.2
\put{$e$}  at 6.87 2.7
\put{$\wedge$} at 6.87 3.2
\put{$\scriptstyle{d+1}$}  at 6.95 3.75
\put{$>$} at 7.7 3.75
\put{$d$}  at 8.2 3.75
\put{$\wedge$} at 8.2 4.3
\put{$c$}  at 8.2 4.9
\put{$\wedge$} at 8.2 5.4
\put{$\vdots$} at 8.2 6.3
\put{$\wedge$} at 8.2 6.7
\put{$b$}  at 8.2 7.4
\put{$\wedge$} at 8.2 8
\put{$a$}  at 8.2 8.6
\putrule from 7.5 9 to 13.5 9               %
\putrule from 6 7.5 to 7.5 7.5              %
\putrule from 10.5 7.5 to 13.5 7.5          
\putrule from 3 6 to 6 6                    %
\putrule from 1.5 3 to 3 3                  %
\putrule from 9 1.5 to 10.5 1.5             %
\putrule from 1.5 0 to 9 0                  %
\putrule from 1.5 0 to 1.5 3            %
\putrule from 3 3 to 3 6                
\putrule from 6 6 to 6 7.5              %
\putrule from 7.5 7.5 to 7.5 9          %
\putrule from 9 0  to  9 1.5            %
\putrule from 10.5 1.5 to 10.5 7.5      %
\putrule from 13.5 7.5 to 13.5 9        %
\putrule from 6.3 0 to 6.3 4.5  
\putrule from 7.5 0 to 7.5 3
\putrule from 6.3 4.5 to 7.5 4.5
\putrule from 7.5 4.5 to 7.5 7.5
\putrule from 7.5 3 to 8.7 3
\putrule from 8.7 3 to 8.7 9 
\endpicture
&
&
\beginpicture
\setcoordinatesystem units <0.5cm,0.5cm>         
\setplotarea x from 0 to 14, y from 0 to 9    
\put{$s_dL=$} at 1 5
\putrule from 7.5 9 to 13.5 9               %
\putrule from 6 7.5 to 7.5 7.5              %
\putrule from 10.5 7.5 to 13.5 7.5          
\putrule from 3 6 to 6 6                    %
\putrule from 1.5 3 to 3 3                  %
\putrule from 9 1.5 to 10.5 1.5             %
\putrule from 1.5 0 to 9 0                  %
\putrule from 1.5 0 to 1.5 3            %
\putrule from 3 3 to 3 6                
\putrule from 6 6 to 6 7.5              %
\putrule from 7.5 7.5 to 7.5 9          %
\putrule from 9 0  to  9 1.5            %
\putrule from 10.5 1.5 to 10.5 7.5      %
\putrule from 13.5 7.5 to 13.5 9        %
\put{$\scriptstyle{d+1}$}  at 4.7 3.8
\put{$d$}  at 4.7 2.8
\putrule from 4.1 2.3 to 4.1 4.3  
\putrule from 5.2 2.3 to 5.2 4.3  
\putrule from 4.1 2.3 to 5.2 2.3
\putrule from 4.1 4.3 to 5.2 4.3
\endpicture
\cr
\cr
\hbox{Case (1)}
&&\hbox{Case (2)} \cr
}
$$
\noindent
Case (1):\   The two boxes in $s_dL$ where the nonstandardness occurs
determine a snake in the shape $\lambda/\mu$ consisting of 
all the boxes above the entry $d$ and all the boxes below the entry $d+1$.   Use
the same notation as in Theorem 5.5 so that $P$
is the standard tableau which is the same as the column reading tableau $C$
except that the entries in the snake are in row reading order.
Then $s_kP$ is nonstandard and there exists a (unique) $v\in S_n$
such that $s_dL = vs_kP$ and $\ell(s_dL) = \ell(v) + \ell(s_kP)$.
The permutation $v$ is given by
$$v = \pmatrix{
i &i+1 &\cdots &k-1 &k &k+1 &k+2 &\cdots &\ell \cr
a &b &\cdots &c &d &d+1 &e &\cdots &f \cr}$$
and $v(x)=L(x)$ for all $x\not\in\{i, i+1, \ldots,\ell\}$.
Thus
$$T_dn_L = T_dT_L n_C = T_vT_kT_Pn_C = T_vT_kn_P.$$ 
Let $\ell(w_{A\cup B})={\ell-i+1\choose 2}$.  By Theorem 5.5
$$T_kn_P = -\sum_{u\le P} (-q)^{\ell(w_{A\cup B})-\ell(w_Aw_B)-\ell(u)}T_u n_C,$$
and so
$$
T_d n_L 
= -T_v\left(\sum_{u\le P} (-q)^{\ell(w_{A\cup B})-\ell(w_Aw_B)-\ell(u)}T_u
n_C\right). 
$$
The permutations $v$ and $u\le P$ are such that $\ell(vu)=\ell(v)+\ell(u)$
and so $T_vT_u=T_{vu}$ for all $u\le P$.  Furthermore, the tableaux
$vuC$, for $u\le P$, are exactly the standard tableaux which are obtained by
permutations of the entries in $L$ which are in the snake.  Thus
$$T_d n_L 
=-\sum_{u\le P} (-q)^{\ell(w_{A\cup B})-\ell(w_Aw_B)-\ell(u)} n_{vuC}.$$

\smallskip\noindent
Case (2):\  Suppose that the boxes containing $d$ and $d+1$ in $L$ are the
boxes containing $k$ and $k+1$ in the column
reading tableau $C$.
$$\beginpicture
\setcoordinatesystem units <0.5cm,0.5cm>         
\setplotarea x from 0 to 14, y from 0 to 9    
\put{$C=$} at 1 5
\putrule from 7.5 9 to 13.5 9               %
\putrule from 6 7.5 to 7.5 7.5              %
\putrule from 10.5 7.5 to 13.5 7.5          
\putrule from 3 6 to 6 6                    %
\putrule from 1.5 3 to 3 3                  %
\putrule from 9 1.5 to 10.5 1.5             %
\putrule from 1.5 0 to 9 0                  %
\putrule from 1.5 0 to 1.5 3            %
\putrule from 3 3 to 3 6                
\putrule from 6 6 to 6 7.5              %
\putrule from 7.5 7.5 to 7.5 9          %
\putrule from 9 0  to  9 1.5            %
\putrule from 10.5 1.5 to 10.5 7.5      %
\putrule from 13.5 7.5 to 13.5 9        %
\put{$k$}  at 4.7 3.8
\put{$\scriptstyle{k+1}$}  at 4.7 2.8
\putrule from 4.1 2.3 to 4.1 4.3  
\putrule from 5.2 2.3 to 5.2 4.3  
\putrule from 4.1 2.3 to 5.2 2.3
\putrule from 4.1 4.3 to 5.2 4.3
\endpicture
$$
Then $Ls_k = s_d L$ and 
$$T_dn_L = T_d T_L n_C = T_L T_k n_C = -q^{-1} T_L n_C = -q^{-1}n_L,$$
by Theorem 5.5 (a).

\smallskip\noindent
(b)  Let $L$ be a standard tableau of shape $\lambda/\mu$ and let $w$ be the
permutation such that $L=wC$, where $C$ is the column reading tableau of shape
$\lambda/\mu$.  By repeatedly using the relations (1.4) we obtain
$$x_in_L = x_i T_w n_C = 
T_w x_{w^{-1}(i)} n_C 
+ \sum_{v<w} T_v p_v(x_1,\ldots, x_n) n_C,$$
where the $p_v(x_1,\ldots, x_n)$ are polynomials in $x_1,\ldots,x_n$
with coefficients in $\ZZ[q,q^{-1}]$.  If $x_i$ acts on $v_C$ by integral
powers of $q$ then we have
$$
x_i n_L = T_w q^{2c(C(w^{-1}(i)))} n_C 
+ \sum_{v<w} T_v p_v(q^{2c(C(1))},\ldots,q^{2c(C(n))}) n_C 
=\sum_{v\le w} b'_v n_{vC}, 
$$
with coefficients $b'_v$ in $\ZZ[q,q^{-1}]$.
\endpf

\bigskip

\section 6.  Induction and restriction

\bigskip

\subsection{Restriction to $H_n$.}

Let $H_n$ be the subalgebra of $\tilde H_n$ generated by
$T_1,\ldots, T_{n-1}$.  The elements $T_w$, $w\in S_n$, form a basis of 
$H_n$.  Since $q$ is not a root of unity the subalgebra $H_n$ of $\tilde
H_n$ is semisimple.  The irreducible representations of $H_n$ are indexed
by the partitions $\nu\vdash n$ and these representations are $q$-analogues
of the irreducible representations of the symmetric group $S_n$.  The following
result describes the decomposition of the restriction to $H_n$ of the
irreducible $\tilde H$-module $\tilde H^{(c,\lambda/\mu)}$.

\thm  Let $\tilde H^{(c,\lambda/\mu)}$ be the irreducible representation
of the affine Hecke algebra $\tilde H_n$ which is defined in Theorem 4.1.
Then
$$\tilde H^{(c,\lambda/\mu)}\big\downarrow^{\tilde H_n}_{H_n}
=\sum_{\nu\vdash n} c_{\mu\nu}^\lambda H^\nu,$$
where $c_{\mu\nu}^\lambda$ is the classical Littlewood-Richardson
coefficient and $H^\nu$ is the irreducible $H_n$-module indexed by the
partition $\nu$. 
\pf
If $\beta=(\beta_1,\ldots,\beta_\ell)$ is a composition of $n$ let
$\gamma_\beta = \gamma_{\beta_1}\times \cdots
\times\gamma_{\beta_\ell}\in S_{\beta_1}\times \cdots\times S_{\beta_\ell}$ where
$\gamma_r = (1,2,\ldots,r)\in S_r$ (in cycle notation).
Let $\chi^{(c,\lambda/\mu)}(T_{\gamma_\beta})$ be the trace of the action of
the element $T_{\gamma_\beta}\in H_n$ on the $\tilde H_n$-module
$\tilde H^{(c,\lambda/\mu)}$.
With notations as in Theorem 4.1
$$\chi^{(c,\lambda/\mu)}(T_{\gamma_{\beta}})
= \sum_{Q\in \cF^{\lambda/\mu}} T_{\gamma_{\beta}}v_Q\big\vert_{v_Q},
$$
and one can copy (without change)
the proof of Theorem 2.20 in [HR2] and obtain
$$\chi^{(c,\lambda/\mu)}(T_{\gamma_{\beta}})
= \sum_{\mu\subseteq\lambda^{(1)}\subseteq\cdots\subseteq
\lambda^{(\ell)} =\lambda} 
\Delta(\lambda^{(1)}/\mu)\Delta(\lambda^{(2)}/\lambda^{(1)})
\cdots \Delta(\lambda^{(\ell)}/\lambda^{(\ell-1)}),$$
where the sum is over all sequences $\mu\subseteq
\lambda^{(1)}\subseteq\cdots\subseteq
\lambda^{(\ell)} =\lambda$ such that $|\lambda^{(i)}/\lambda^{(i-1)}|=\beta_i$
and the factor $\Delta(\lambda^{(i)}/\lambda^{(i-1)})$ is given by
$$\Delta(\lambda/\mu)
= \cases{
\displaystyle{
(q-q^{-1})^{cc-1}\prod_{bs\in CC} q^{c(bs)-1}(-q^{-1})^{r(bs)-1},
}
&if $\lambda/\mu$ is a border strip, \cr
0, &otherwise. \cr}$$
In the formula for $\Delta(\lambda/\mu)$:
a border strip is a skew shape with at
most one box in each diagonal, $CC$ is the set of connected components
of $\lambda/\mu$, $cc$ is the number of connected components of $\lambda/\mu$,
$r(bs)$ is the number of rows of $bs$, and $c(bs)$ is the number of columns of
$bs$.

Let
$s_\lambda$ denote the Schur function (see [Mac]) and define
$$q_r = \sum_{m=1}^r (-q^{-1})^{r-m}q^{m-1}s_{(m1^{r-m})}.$$
By Proposition 6.11(a) in [HR1],
$$q_r s_\mu = \sum_{\lambda} \Delta(\lambda/\mu) s_\lambda.$$
Letting $q_\beta = q_{\beta_1}\cdots q_{\beta_\ell}$ one can inductively
apply this formula to obtain
$$q_\beta s_\mu = \sum_{\lambda}
\chi^{(c,\lambda/\mu)}(T_{\gamma_{\beta}}) s_\lambda.$$
Thus, with notations as in [Mac] Chapt. I,
$$\eqalign{
\chi^{(c,\lambda/\mu)}(T_{\gamma_\beta})
&=\langle q_\beta s_\mu, s_\lambda\rangle \cr
&= \langle q_\beta , s_{\lambda/\mu}\rangle,
\qquad\qquad\qquad\ \hbox{by [Mac] I (5.1),} \cr 
&= \sum_{\nu} c_{\mu\nu}^\lambda \langle q_\beta, s_\nu\rangle,
\qquad\qquad\ \hbox{by [Mac] I (5.3),} \cr
&= \sum_{\nu} c_{\mu\nu}^\lambda \chi^\nu(T_{\gamma_\beta}),  
\qquad\qquad\,\hbox{by [Ra1] Theorem 4.14,} \cr}$$
where $\chi^\nu(T_{\gamma_\beta})$ denotes the irreducible character
of $H_n$ evaluated at the element $T_{\gamma_\beta}$.
The result follows since, by [Ra1] Theorem 5.1,
the characters of $H_n$ are determined by their values on the elements
$T_{\gamma_\beta}$. 
\endpf

\noindent
Classically, the Littlewood-Richardson
coefficients describe
\smallskip
\item{(1)} The decomposition of the tensor product of two irreducible
polynomial representations of $GL_n(\CC)$, and
\smallskip
\item{(2)} The decomposition of an irreducible representation of
$S_k\times S_\ell$ when it is induced to $S_{k+\ell}$.
\smallskip\noindent
Theorem 6.2 gives an exciting {\it new} way of interpreting these coefficients.
They describe
\item{(3)}  The decomposition of an irreducible representation $\tilde
H_n$ when it is restricted to the subalgebra $H_n$.
\endprop

\subsection{Induction from ``Young subalgebras''.}

Let $k$ and $\ell$ be such that $k+\ell=n$.  Let
$$
\eqalign{
\tilde H_k &= \hbox{ the subalgebra of $\tilde H_n$ generated by
$T_i, 1\le i\le k-1$ and $x_i$, $1\le i\le k$,} \cr
\tilde H_\ell &= \hbox{ the
subalgebra of $\tilde H_n$ generated by $T_i$, $k+1\le i\le n-1$, and
$x_i$, $k+1\le i\le n$.} \cr
}$$
In this way $\tilde H_k\otimes \tilde H_\ell$ is 
naturally a subalgebra of $\tilde H_n$.

Let $(a,\theta)$  be a placed skew shape with $k$ boxes and let $(b,\phi)$ be a
placed skew shape with $\ell$ boxes.  Number the boxes of $\theta$ with
$1,\ldots, k$ (as in Section 2, along diagonals from southwest to northeast) 
and number the boxes of $\phi$ with
$k+1,\ldots, n$, in order to match the imbeddings of $\tilde H_k$ and
$\tilde H_\ell$ in $\tilde H_n$.  Let $\tilde H^{(a,\theta)}$ and
$\tilde H^{(b,\phi)}$ be the corresponding representations of
$\tilde H_k$ and $\tilde H_\ell$ as defined in Theorem 4.1.

Let $\theta *_v\phi$ (resp. $\theta*_h \phi$)
be the skew shape obtained by placing
$\theta$ and $\phi$ adjacent to each other in such a way that ${\rm box}_{(k+1)}$
of $\phi$ is immediately above (resp. to the left of) 
${\rm box}_k$ of $\theta$.  Let $a\otimes b$ be the content function
given by
$$(a\otimes b)({\rm box}_i) = \cases{
a({\rm box}_i), &if $1\le i\le k$, \cr
b({\rm box}_i), &if $k+1\le i\le \ell$, \cr}
$$
\thm 
With notations as above,
$$\Ind_{\tilde H_k\otimes \tilde H_\ell}^{\tilde H_n}
(\tilde H^{(a,\theta)}\otimes \tilde H^{(b,\phi)})
=\tilde H^{(a\otimes b,\theta*_v\phi)} + \tilde H^{(a\otimes b, \theta*_h\phi)}$$
in the Grothendieck ring of finite dimensional representations of $\tilde
H_n$.
\pf  Let $S_n/(S_k\times S_\ell)$ be the set of minimal length coset
representatives of the cosets of $S_k\times S_\ell$ in $S_n$.  
The module $M=\Ind_{\tilde H_k\otimes \tilde H_\ell}^{\tilde H_n}
(\tilde H^{(a,\theta)}\otimes \tilde H^{(b,\phi)})
$ has basis
$$T_w(v_L\otimes v_Q)
\qquad\hbox{
where $w\in S_n/(S_k\times S_\ell)$, $L\in \cF^{\theta}$
and $Q\in \cF^{\phi}$.}$$ 
By repeatedly applying the relations (1.4) we obtain
$$x_i(T_w(v_L\otimes v_Q))
=T_w x_{w^{-1}(i)}(v_L\otimes v_Q) + \sum_{u<w} b_u T_u(v_L\otimes v_Q),$$
for some constants $b_u\in \CC$.  From this we can see that
the action of $x_i$ on $M$ is an upper triangular matrix
with eigenvalues $q^{2c(P(i))}$ where $P$ runs over
the standard tableaux of shapes
$\theta*_v\phi$ and $\theta*_h\phi$.  It follows that $M$ has 
$${|S_n|\Card(\cF^{\theta})\Card(\cF^{\phi})\over
|S_k\times S_\ell|}$$
distinct weights.  Since this number is exactly the dimension of
$M$, it follows that every generalized weight space of $M$ is one dimensional
and thus that $M$ is calibrated.
By Theorem 4.1, all irreducible calibrated representations are of the form
$\tilde H^{(c,\lambda/\mu)}$ for some placed skew shape $(c,\lambda/\mu)$ and by
Lemma 2.2 this placed skew shape is completely determined
by any one of the weights of the module $\tilde H^{(c,\lambda/\mu)}$.  Thus,
our analysis of the weights of $M$ implies that
both $\tilde H^{(a\otimes b,\theta*_v\phi)}$ and $\tilde H^{(a\otimes
b,\theta*_h\phi)}$ are composition factors of $M$.  The result follows since
$$\dim(\tilde H^{(a\otimes b,\theta*_v\phi)}) 
+ \dim(\tilde H^{(a\otimes b,\theta*_h\phi)}) =
\dim(M).\qquad\hbox{\qed}$$
\medskip

A {\it ribbon} is a skew shape which has at most one box in each
diagonal.

\cor Let $c$ be the content function given by
$c({\rm box}_i) = i-1$, for $1\le i\le n$.
Let $t=(t_1,\ldots,t_n)=(1,q^2,\ldots,q^{2(n-1)})$ and let
$\CC v_t$ be the one dimensional module for
$\CC[X]=\CC[x_1^{\pm1},\ldots,x_n^{\pm1}]$ given by
$$x_iv_t = t_iv_t,\qquad\hbox{for all $1\le i\le n$}.$$
In the Grothendieck ring of finite dimensional $\tilde H_n$-modules
$$\Ind_{\CC[X]}^{\tilde H_n}(\CC v_t) =
\sum_{\lambda/\mu} \tilde H^{(c,\lambda/\mu)},$$
where the sum is over all connected ribbons $\lambda/\mu$ with $n$ boxes.
\pf
Since $\CC[X] = \tilde H_1\otimes\tilde H_1\otimes\cdots\otimes \tilde H_1
\subseteq \tilde H_n$ this result can be obtained by repeatedly applying
Theorem 6.2.
\endpf

\medskip\noindent
Theorem 6.2 and Corollary 6.3 are $\tilde H_n$-module realizations of
the Schur function identities in [Mac] I \S 5 Ex 21 (a),(b).
The module $\Ind_{\CC[X]}^{\tilde H_n}(\CC v_t)$ is a
principal series module for $\tilde H_n$ (see [Ka]).  The identity in 
Corollary 6.3 describes the composition series of this principal series module.
Using the methods of [Ra3] and [Ra4, (1.2) Ex. 2] one can obtain a
generalization of this identity (and thus of [Mac] I \S 5 Ex 21 (b)) which holds
for affine Hecke algebras of arbitrary Lie type.
\endprop






\vfill\eject

\centerline{\smallcaps References}
\bigskip

\medskip
\item{[AK]} {\smallcaps S.\ Ariki and K.\ Koike}, {\it A Hecke algebra of
$(\ZZ/r\ZZ)\wr S_n$ and construction of its irreducible representations},
Adv.\ in Math.\ {\bf 106} (1994), 216--243.

\medskip
\item{[BW]} {\smallcaps A.\ Bj\"orner and M.\ Wachs}, {\it
Generalized quotients in Coxeter groups}, Trans. Amer. Math. Soc. {\bf 308}
(1988), 1--37. 

\medskip
\item{[Bou]} {\smallcaps N.\ Bourbaki}, 
{\it Groupes et alg\`ebres de Lie, Chapitres 4,5 et 6}, 
Elements de Math\'e\-matique, Hermann, Paris 1968.

\medskip
\item{[Ch]} {\smallcaps I.\ Cherednik}, {\it A new
interpretation of Gel'fand-Tzetlin bases}, Duke Math. J. {\bf 54} (1987),
563--577.

\medskip
\item{[Cr]} {\smallcaps J.\ Crisp}, Ph.D. Thesis, University of Sydney, 1997.

\medskip
\item{[D]} {\smallcaps P.\ Diaconis}, Lecture at the {\sl Workshop on
Representation Theory and Symmetric Functions}, Mathematical Sciences Research
Institute, Berkeley, April 1997.

\medskip
\item{[Fo]} {\smallcaps S.\ Fomin}, personal communication, 1997.

\medskip
\item{[Fu]} {\smallcaps W.\ Fulton}, {\sl Young tableaux: With
applications to representation theory and geometry}, London Mathematical Society
Student Texts {\bf 35}, Cambridge University Press, Cambridge, 1997.

\medskip
\item{[GL]} {\smallcaps M.\ Geck and S.\ Lambropoulou}, {\it Markov
traces and knot invariants related to Iwahori-Hecke algebras of type $B$}, J.
Reine Angew. Math. {\bf 482} (1997), 191--213.

\medskip
\item{[GW]} {\smallcaps A.\ Garsia, M.\ Wachs}, {\it Combinatorial aspects of
skew representations of the symmetric group}, J. Combin. Theory Ser. A {\bf 50}
(1989), 47--81. 

\medskip
\item{[HR1]} {\smallcaps T.\ Halverson and A.\ Ram},
{\it Characters of algebras containing a Jones basic construction: the 
Temperley-Lieb, Okada, Brauer, and Birman-Wenzl algebras}, Adv. Math. {\bf 116}
(1995), 263--321.

\medskip
\item{[HR2]} {\smallcaps T.\ Halverson and A.\ Ram},
{\it Murnaghan-Nakayama rules for characters of Iwahori-Hecke algebras of
classical type}, Trans. Amer. Math. Soc. {\bf 348} (1996), 3967--3995. 

\medskip
\item{[Ho]} {\smallcaps P.N.\ Hoefsmit}, 
{\it Representations of Hecke algebras of finite groups with $BN$-pairs 
of classical type}, 
Ph.D.\ Thesis, University of British Columbia, 1974.

\medskip
\item{[JK]}  {\smallcaps G.\ James and A.\  Kerber}, {\sl The
representation theory of the symmetric group}, Encyclopedia of Mathematics and
its Applications {\bf 16}, Addison-Wesley Publishing Co., Reading, Mass., 1981.

\medskip
\item{[Jo]} {\smallcaps V.F.R.\ Jones}, {\it A quotient of the affine Hecke
algebra in the Brauer algebra}, Enseign. Math. (2) {\bf 40} (1994), 313--344. 

\medskip
\item{[Ka]} {\smallcaps S-i.\ Kato},
{\it Irreducibility of principal series representations for Hecke
algebras of affine type}, 
J. Fac. Sci. Univ. Tokyo Sect. IA Math. {\bf 28} (1981), 929--943.

\medskip
\item{[Lu]} {\smallcaps G.\ Lusztig},
{\it Affine Hecke algebras and their graded version}, J. Amer. Math.
Soc. {\bf 2} (1989), 599--635.

\medskip
\item{[Mac]} {\smallcaps I.G.\ Macdonald}, 
{\sl Symmetric functions and Hall polynomials}, Second edition, Oxford
University Press, New York, 1995. 

\medskip
\item{[Mac2]} {\smallcaps I.G.\ Macdonald},
{\sl Spherical functions on a group of $p$-adic type}, Ramanujan Insititute
for Advanced Study, University of Madras, Madras, India, 1971.

\medskip
\item{[Op]}  {\smallcaps E.\ Opdam},
{\it Harmonic analysis for certain representations of graded Hecke algebras},
Acta Math. {\bf 175} (1995), 75--121.

\medskip
\item{[Ra1]} {\smallcaps A.\ Ram}, {\it A Frobenius formula for the
characters of the Hecke algebras}, Invent. Math. {\bf 106} (1991),
461--488.

\medskip
\item{[Ra2]} {\smallcaps A.\ Ram}, 
{\it Seminormal representations of Weyl groups and Iwahori-Hecke algebras}, 
Proc.\ London Math.\ Soc.\ (3) {\bf 75} (1997), 99-133.

\medskip
\item{[Ra3]} {\smallcaps A.\ Ram}, 
{\it Calibrated representations of affine Hecke algebras}, preprint 1998.

\medskip
\item{[Ra4]} {\smallcaps A. \ Ram}, 
{\it Standard Young tableaux for finite root systems}, preprint 1998.

\medskip
\item{[Ra5]} {\smallcaps A.\ Ram}, 
{\it Irreducible representations
of rank two affine Hecke algebras}, preprint 1998.

\medskip
\item{[RR1]} {\smallcaps A.\ Ram and J. Ramagge}, 
{\it  Jucys-Murphy elements come from affine Hecke algebras}, in preparation.

\medskip
\item{[RR2]} {\smallcaps A.\ Ram and J. Ramagge}, {\it Calibrated
representations and the $q$-Springer correspondence}, in preparation.

\medskip
\item{[Sg]} {\smallcaps B.\ Sagan}, {\sl The symmetric group,
Representations, combinatorial algorithms, and symmetric functions}, The
Wadsworth \& Brooks/Cole Mathematics Series, Wadsworth \& Brooks/Cole Advanced
Books \& Software, Pacific Grove, CA, 1991.

\medskip
\item{[Wz]}  {\smallcaps H. Wenzl}, 
{\it Hecke algebras of type $A\sb n$ and subfactors}, Invent. Math. {\bf 92}
(1988), 349--383.

\medskip
\item{[Y]}  {\smallcaps A.\ Young}, {\it On quantitative substitutional analysis}
(sixth and eighth papers),  Proc. London Math. Soc. (2) {\bf 34} (1931), 196--230
and {\bf 37} (1934), 441--495.
\vfill\eject
\end